\documentclass[11pt]{article}
\pdfoutput=1 
\usepackage[utf8x]{inputenc}
\usepackage[english]{babel}
\usepackage{graphicx}           
\usepackage{url}
\usepackage{eurosym}
\usepackage{makeidx}
\usepackage[plain]{algorithm}
\usepackage{algorithmic} 
\usepackage{color} 
\usepackage[margin=3.3cm]{geometry}
\usepackage{setspace}

\usepackage{amsmath,amsfonts,amssymb,latexsym,stmaryrd}

\let\oldmarginpar\marginpar
\renewcommand\marginpar[1]{\-\oldmarginpar[\raggedleft\footnotesize #1]%
{\raggedright\footnotesize #1}}

\begin{document}

\typeout{>>> macros latex 17 septembre 2002}

\newsavebox{\fmbox}
\newenvironment{fabbox}
     {\noindent\begin{lrbox}{\fmbox}\begin{minipage}{0.983\textwidth}}
     {\end{minipage}\end{lrbox}\fbox{\usebox{\fmbox}}}

\newenvironment{mat}
     {\left(\begin{smallmatrix}}
     {\end{smallmatrix}\right)}

\newcommand{\mmin}{{\textrm{min}}}
\newcommand{\mmax}{{\textrm{max}}}
\newcommand{\var}{{\textrm{var}}}

\newcommand{\V}        {\mathbb V}
\newcommand{\M}        {\mathbb M}
\newcommand{\N}        {\mathbb N}
\newcommand{\D}        {\mathbb D}
\newcommand{\Z}        {\mathbb Z}
\newcommand{\C}        {\mathbb C}
\newcommand{\T}        {\mathbb T}
\newcommand{\E}        {\mathbb E}
\newcommand{\X}        {\mathbb X}
\newcommand{\Y}        {\mathbb Y}
\newcommand{\R}      {\mathbb R}
\newcommand{\Q}      {\mathbb Q}
\newcommand{\Para}    {\P}
\renewcommand{\P}      {\mathbb P}

\newcommand{\rfr}[1]    {\stackrel{\circ}{#1}}
\newcommand{\equiva}    {\displaystyle\mathop{\simeq}}
\newcommand{\equivb}    {\displaystyle\mathop{\sim}}
\newcommand{\simiid}     {\stackrel{\textrm{iid}}{\sim}}

\newcommand{\gtlt}     {\displaystyle\mathop{\gtrless}}

\newcommand{\Limsup}    {\mathop{\overline{\textrm lim}}}
\newcommand{\Liminf}    {\mathop{\underline{\textrm lim}}}
\newcommand{\osc}       {\mathop{\hbox{\textrm osc}}}
\newcommand{\ccv}       {\mathop{\;{\longrightarrow}\;}}
\newcommand{\cv}        {\mathop{\;{\rightarrow}\;}}
\newcommand{\cvweak}    {\mathop{\;{\rightharpoonup}\;}}
\newcommand{\versup}    {\mathop{\;{\nearrow}\;}}
\newcommand{\versdown}  {\mathop{\;{\searrow}\;}}
\newcommand{\vvers}     {\mathop{\;{\longrightarrow}\;}}
\newcommand{\cvetroite} {\mathop{\;{\Longrightarrow}\;}}
\newcommand{\cvlaw}     {\mathop{\;{\Longrightarrow}\;}}
\newcommand{\argmax}    {\hbox{{\textrm Arg}}\max}
\newcommand{\argmin}    {\hbox{{\textrm Arg}}\min}
\newcommand{\esssup}    {\hbox{{\textrm ess}}\sup}
\newcommand{\essinf}    {\hbox{{\textrm ess}}\inf}
\renewcommand{\div}     {\hbox{{\textrm div}}}
\newcommand{\rot}       {\hbox{{\textrm rot}}}
\newcommand{\supp}      {\hbox{{\textrm supp}}}
\newcommand{\cov}       {\hbox{{\textrm cov}}}
\newcommand{\indep}      {\perp\!\!\!\!\perp}
\newcommand{\abs}    [1] {\left| #1 \right|}
\newcommand{\norm}   [1] {\left\Vert #1 \right\Vert}
\newcommand{\normtv}   [1] {\left\Vert #1 \right\Vert_{\textrm{\tiny TV}}}
\newcommand{\Norm}   [1] {\Vert #1 \Vert}
\newcommand{\Normtv}   [1] {\Vert #1 \Vert_{\textrm{\tiny TV}}}
\newcommand{\scrochet}[1] {{\tiny\langle} #1 {\scriptsize\rangle}}
\newcommand{\crochet}[1] {\langle #1 \rangle}
\newcommand{\bigcrochet}[1] {\big\langle #1 \big\rangle}
\newcommand{\Bigcrochet}[1] {\Big\langle #1 \Big\rangle}
\newcommand{\bcrochet}[1] {\boldsymbol{\langle} #1 \boldsymbol{\rangle}}
\newcommand{\bigbcrochet}[1] {\boldsymbol{\big\langle} #1 \boldsymbol{\big\rangle}}
\newcommand{\Bigbcrochet}[1] {\boldsymbol{\Big\langle} #1 \boldsymbol{\Big\rangle}}
\newcommand{\Crochet}[1] {\left\langle #1 \right\rangle}

\newcommand{\espc}   [3] {E_{#1}\left(\left. #2 \right| #3 \right)}
\newcommand{\trace}     {\hbox{{\textrm trace}}}
\newcommand{\diag}      {\hbox{{\textrm diag}}}
\newcommand{\tbullet}   {$\bullet$}
\newcommand{\ot}        {\leftarrow}
\newcommand{\carre}     {\hfill$\Box$}
\newcommand{\carreb}    {\hfill\rule{0.25cm}{0.25cm}}
\newcommand{\demi}{{{\textstyle\frac{1}{2}}}}
\newcommand{\quart}{{{\textstyle\frac{1}{4}}}}

\newcommand{\ttA}  {{\texttt A}}
\newcommand{\ttB}  {{\texttt B}}
\newcommand{\ttC}  {{\texttt C}}
\newcommand{\ttD}  {{\texttt D}}
\newcommand{\ttE}  {{\texttt E}}
\newcommand{\ttF}  {{\texttt F}}
\newcommand{\ttG}  {{\texttt G}}
\newcommand{\ttH}  {{\texttt H}}
\newcommand{\ttI}  {{\texttt I}}
\newcommand{\ttJ}  {{\texttt J}}
\newcommand{\ttK}  {{\texttt K}}
\newcommand{\ttL}  {{\texttt L}}
\newcommand{\ttM}  {{\texttt M}}
\newcommand{\ttN}  {{\texttt N}}
\newcommand{\ttO}  {{\texttt O}}
\newcommand{\ttP}  {{\texttt P}}
\newcommand{\ttQ}  {{\texttt Q}}
\newcommand{\ttR}  {{\texttt R}}
\newcommand{\ttS}  {{\texttt S}}
\newcommand{\ttT}  {{\texttt T}}
\newcommand{\ttU}  {{\texttt U}}
\newcommand{\ttV}  {{\texttt V}}
\newcommand{\ttW}  {{\texttt W}}
\newcommand{\ttX}  {{\texttt X}}
\newcommand{\ttY}  {{\texttt Y}}
\newcommand{\ttZ}  {{\texttt Z}}
\newcommand{\tta}  {{\texttt a}}
\newcommand{\ttb}  {{\texttt b}}
\newcommand{\ttc}  {{\texttt c}}
\newcommand{\ttd}  {{\texttt d}}
\newcommand{\tte}  {{\texttt e}}
\newcommand{\ttf}  {{\texttt f}}
\newcommand{\ttg}  {{\texttt g}}
\newcommand{\tth}  {{\texttt h}}
\newcommand{\tti}  {{\texttt i}}
\newcommand{\ttj}  {{\texttt j}}
\newcommand{\ttk}  {{\texttt k}}
\newcommand{\ttl}  {{\texttt l}}
\newcommand{\ttm}  {{\texttt m}}
\newcommand{\ttn}  {{\texttt n}}
\newcommand{\tto}  {{\texttt o}}
\newcommand{\ttp}  {{\texttt p}}
\newcommand{\ttq}  {{\texttt q}}
\newcommand{\ttr}  {{\texttt r}}
\newcommand{\tts}  {{\texttt s}}
\newcommand{\ttt}  {{\texttt t}}
\newcommand{\ttu}  {{\texttt u}}
\newcommand{\ttv}  {{\texttt v}}
\newcommand{\ttw}  {{\texttt w}}
\newcommand{\ttx}  {{\texttt x}}
\newcommand{\tty}  {{\texttt y}}
\newcommand{\ttz}  {{\texttt z}}

\newcommand{\sfA}  {{\mathsf A}}
\newcommand{\sfB}  {{\mathsf B}}
\newcommand{\sfC}  {{\mathsf C}}
\newcommand{\sfD}  {{\mathsf D}}
\newcommand{\sfE}  {{\mathsf E}}
\newcommand{\sfF}  {{\mathsf F}}
\newcommand{\sfG}  {{\mathsf G}}
\newcommand{\sfH}  {{\mathsf H}}
\newcommand{\sfI}  {{\mathsf I}}
\newcommand{\sfJ}  {{\mathsf J}}
\newcommand{\sfK}  {{\mathsf K}}
\newcommand{\sfL}  {{\mathsf L}}
\newcommand{\sfM}  {{\mathsf M}}
\newcommand{\sfN}  {{\mathsf N}}
\newcommand{\sfO}  {{\mathsf O}}
\newcommand{\sfP}  {{\mathsf P}}
\newcommand{\sfQ}  {{\mathsf Q}}
\newcommand{\sfR}  {{\mathsf R}}
\newcommand{\sfS}  {{\mathsf S}}
\newcommand{\sfT}  {{\mathsf T}}
\newcommand{\sfU}  {{\mathsf U}}
\newcommand{\sfV}  {{\mathsf V}}
\newcommand{\sfW}  {{\mathsf W}}
\newcommand{\sfX}  {{\mathsf X}}
\newcommand{\sfY}  {{\mathsf Y}}
\newcommand{\sfZ}  {{\mathsf Z}}
\newcommand{\sfa}  {{\mathsf a}}
\newcommand{\sfb}  {{\mathsf b}}
\newcommand{\sfc}  {{\mathsf c}}
\newcommand{\sfd}  {{\mathsf d}}
\newcommand{\sfe}  {{\mathsf e}}
\newcommand{\sff}  {{\mathsf f}}
\newcommand{\sfg}  {{\mathsf g}}
\newcommand{\sfh}  {{\mathsf h}}
\newcommand{\sfi}  {{\mathsf i}}
\newcommand{\sfj}  {{\mathsf j}}
\newcommand{\sfk}  {{\mathsf k}}
\newcommand{\sfl}  {{\mathsf l}}
\newcommand{\sfm}  {{\mathsf m}}
\newcommand{\sfn}  {{\mathsf n}}
\newcommand{\sfo}  {{\mathsf o}}
\newcommand{\sfp}  {{\mathsf p}}
\newcommand{\sfq}  {{\mathsf q}}
\newcommand{\sfr}  {{\mathsf r}}
\newcommand{\sfs}  {{\mathsf s}}
\newcommand{\sft}  {{\mathsf t}}
\newcommand{\sfu}  {{\mathsf u}}
\newcommand{\sfv}  {{\mathsf v}}
\newcommand{\sfw}  {{\mathsf w}}
\newcommand{\sfx}  {{\mathsf x}}
\newcommand{\sfy}  {{\mathsf y}}
\newcommand{\sfz}  {{\mathsf z}}

\newcommand{\bfA}  {{\mathbf A}}
\newcommand{\bfB}  {{\mathbf B}}
\newcommand{\bfC}  {{\mathbf C}}
\newcommand{\bfD}  {{\mathbf D}}
\newcommand{\bfE}  {{\mathbf E}}
\newcommand{\bfF}  {{\mathbf F}}
\newcommand{\bfG}  {{\mathbf G}}
\newcommand{\bfH}  {{\mathbf H}}
\newcommand{\bfI}  {{\mathbf I}}
\newcommand{\bfJ}  {{\mathbf J}}
\newcommand{\bfK}  {{\mathbf K}}
\newcommand{\bfL}  {{\mathbf L}}
\newcommand{\bfM}  {{\mathbf M}}
\newcommand{\bfN}  {{\mathbf N}}
\newcommand{\bfO}  {{\mathbf O}}
\newcommand{\bfP}  {{\mathbf P}}
\newcommand{\bfQ}  {{\mathbf Q}}
\newcommand{\bfR}  {{\mathbf R}}
\newcommand{\bfS}  {{\mathbf S}}
\newcommand{\bfT}  {{\mathbf T}}
\newcommand{\bfU}  {{\mathbf U}}
\newcommand{\bfV}  {{\mathbf V}}
\newcommand{\bfW}  {{\mathbf W}}
\newcommand{\bfX}  {{\mathbf X}}
\newcommand{\bfY}  {{\mathbf Y}}
\newcommand{\bfZ}  {{\mathbf Z}}
\newcommand{\bfa}  {{\mathbf a}}
\newcommand{\bfb}  {{\mathbf b}}
\newcommand{\bfc}  {{\mathbf c}}
\newcommand{\bfd}  {{\mathbf d}}
\newcommand{\bfe}  {{\mathbf e}}
\newcommand{\bff}  {{\mathbf f}}
\newcommand{\bfg}  {{\mathbf g}}
\newcommand{\bfh}  {{\mathbf h}}
\newcommand{\bfi}  {{\mathbf i}}
\newcommand{\bfj}  {{\mathbf j}}
\newcommand{\bfk}  {{\mathbf k}}
\newcommand{\bfl}  {{\mathbf l}}
\newcommand{\bfm}  {{\mathbf m}}
\newcommand{\bfn}  {{\mathbf n}}
\newcommand{\bfo}  {{\mathbf o}}
\newcommand{\bfp}  {{\mathbf p}}
\newcommand{\bfq}  {{\mathbf q}}
\newcommand{\bfr}  {{\mathbf r}}
\newcommand{\bfs}  {{\mathbf s}}
\newcommand{\bft}  {{\mathbf t}}
\newcommand{\bfu}  {{\mathbf u}}
\newcommand{\bfv}  {{\mathbf v}}
\newcommand{\bfw}  {{\mathbf w}}
\newcommand{\bfx}  {{\mathbf x}}
\newcommand{\bfy}  {{\mathbf y}}
\newcommand{\bfz}  {{\mathbf z}}

\newcommand{\bbA}  {{\mathbb A}}
\newcommand{\bbB}  {{\mathbb B}}
\newcommand{\bbC}  {{\mathbb C}}
\newcommand{\bbD}  {{\mathbb D}}
\newcommand{\bbE}  {{\mathbb E}}
\newcommand{\bbF}  {{\mathbb F}}
\newcommand{\bbG}  {{\mathbb G}}
\newcommand{\bbH}  {{\mathbb H}}
\newcommand{\bbI}  {{\mathbb I}}
\newcommand{\bbJ}  {{\mathbb J}}
\newcommand{\bbK}  {{\mathbb K}}
\newcommand{\bbL}  {{\mathbb L}}
\newcommand{\bbM}  {{\mathbb M}}
\newcommand{\bbN}  {{\mathbb N}}
\newcommand{\bbO}  {{\mathbb O}}
\newcommand{\bbP}  {{\mathbb P}}
\newcommand{\bbQ}  {{\mathbb Q}}
\newcommand{\bbR}  {{\mathbb R}}
\newcommand{\bbS}  {{\mathbb S}}
\newcommand{\bbT}  {{\mathbb T}}
\newcommand{\bbU}  {{\mathbb U}}
\newcommand{\bbV}  {{\mathbb V}}
\newcommand{\bbW}  {{\mathbb W}}
\newcommand{\bbX}  {{\mathbb X}}
\newcommand{\bbY}  {{\mathbb Y}}
\newcommand{\bbZ}  {{\mathbb Z}}

\newcommand{\AAA}  {{\mathcal A}}
\newcommand{\BB}   {{\mathcal B}}
\newcommand{\CC}   {{\mathcal C}}
\newcommand{\DD}   {{\mathcal D}}
\newcommand{\EE}   {{\mathcal E}}
\newcommand{\FF}   {{\mathcal F}}
\newcommand{\GG}   {{\mathcal G}}
\newcommand{\HH}   {{\mathcal H}}
\newcommand{\II}   {{\mathcal I}}
\newcommand{\JJ}   {{\mathcal J}}
\newcommand{\KK}   {{\mathcal K}}
\newcommand{\LL}   {{\mathcal L}}
\newcommand{\NN}   {{\mathcal N}}
\newcommand{\MM}   {{\mathcal M}}
\newcommand{\OO}   {{\mathcal O}}
\newcommand{\PP}   {{\mathcal P}}
\newcommand{\QQ}   {{\mathcal Q}}
\newcommand{\RR}   {{\mathcal R}}
\renewcommand{\SS}   {{\mathcal S}}
\newcommand{\SSS}   {{\mathcal S}}
\newcommand{\TT}   {{\mathcal T}}
\newcommand{\UU}   {{\mathcal U}}
\newcommand{\VV}   {{\mathcal V}}
\newcommand{\WW}   {{\mathcal W}}
\newcommand{\XX}   {{\mathcal X}}
\newcommand{\YY}   {{\mathcal Y}}
\newcommand{\ZZ}   {{\mathcal Z}}


\newcommand{\rmd}   {{{\textrm{\upshape d}}}}

\newcommand{\dontforget}[1]
{{\mbox{}\\\noindent\rule{1cm}{2mm}\hfill  #1 \hfill\rule{1cm}{2mm}}\typeout{---------- #1 ------------}}
\newcommand{\rond}[1]     {\stackrel{\circ}{#1}}
\newcommand{\indic}{{\mathrm\mathbf1}}
\renewcommand{\epsilon}{\varepsilon}
\newcommand{\marginal}[1]{%
        \leavevmode\marginpar{\tiny\raggedright#1\par}}
\newcommand{\warning}{\setlength{\unitlength}{1cm}
  \begin{picture}(0.6,0.5)(0,0)
  \put(0.25,0.15){\makebox(0,0){{\huge$\bigtriangleup$}}}
  \put(0.25,0.16){\makebox(0,0){{\small\sf !}}}
  \end{picture}}

\def\bm#1{%
  \mathchoice%
       {\setbox1=\hbox{$#1$}\dobm}
       {\setbox1=\hbox{$#1$}\dobm}
       {\setbox1=\hbox{\scriptsize$#1$}\dobm}
       {\setbox1=\hbox{\tiny$#1$}\dobm}}
\def\dobm{
    \copy1\kern-\wd1\kern0.05ex\copy1\kern-\wd1\kern0.05ex\box1}

\renewcommand{\theenumi}{\roman{enumi}}
\renewcommand{\labelenumi}{{\textrm{\rm({\it\theenumi}\/)}}}
\renewcommand{\theenumii}{\alph{enumii}}
\renewcommand{\labelenumii}{{\textit\theenumii.}}

\newcommand{\fenumi}  {\textrm{\rm({\textit{i}}\/)}}
\newcommand{\fenumii} {\textrm{\rm({\textit{ii}}\/)}}
\newcommand{\fenumiii}{\textrm{\rm({\textit{iii}}\/)}}
\newcommand{\fenumiv} {\textrm{\rm({\textit{iv}}\/)}}
\newcommand{\fenumv}  {\textrm{\rm({\textit{v}}\/)}}
\newcommand{\fenum}[1]{\textrm{\rm({\textit{#1}}\/)}}

\newcommand
      {\sysdys}
      {{\sf S\kern-.15em\raise.3ex\hbox{Y}\kern-.15em
            SD\kern-.15em\raise.3ex\hbox{Y}\kern-.15emS}}

%
\newcommand{\eqdef}     {\stackrel{{\textrm{\rm\tiny def}}}{=}}

\newtheorem{theorem}      {Theorem}[section]
\newtheorem{theorem*}     {theorem}
\newtheorem{proposition}  [theorem]{Proposition}
\newtheorem{definition}   [theorem]{Definition}
\newtheorem{lemma}        [theorem]{Lemma}
\newtheorem{corollary}    [theorem]{Corollary}
\newtheorem{result}       [theorem]{Result}
\newtheorem{hypothesis}   [theorem]{Hypothesis}
\newtheorem{hypotheses}   [theorem]{Hypotheses}
\newtheorem{remark}       [theorem]{Remark}
\newtheorem{remarks}      [theorem]{Remarks}
\newtheorem{property}     [theorem]{Property}
\newtheorem{properties}   [theorem]{Properties}
\newtheorem{example}      [theorem]{Example}
\newtheorem{examples}     [theorem]{Examples}
%
%
\newcommand{\proof}        {\paragraph{Proof}}

\newcommand{\eqas}      {\stackrel{\hbox{{\scriptsize a.s.}}}{=}}
\newcommand{\eqlaw}     {\stackrel{\hbox{{\scriptsize law}}}{=}}
\newcommand{\law}       {{\rm law}}
\newcommand{\limas}     {\mathop{\hbox{\rm lim--a.s.}}}
\newcommand{\rank}       {\hbox{rank}}
\newcommand{\sign}      {\hbox{{\rm sign}}}


\newcommand{\lambdad}{\lambda^{\textrm{\tiny d}}}
\newcommand{\lambdab}{\lambda^{\textrm{\tiny b}}}
\newcommand{\lambdac}{\lambda^{\textrm{\tiny c}}}
\newcommand{\lambdadm}
           {\lambda^{\textrm{\tiny d}}_{\textrm{\tiny max}}}
\newcommand{\lambdabm}
           {\lambda^{\textrm{\tiny b}}_{\textrm{\tiny max}}}
\newcommand{\lambdacm}
           {\lambda^{\textrm{\tiny c}}_{\textrm{\tiny max}}}
\newcommand{\umax}
           {u_{\textrm{\tiny max}}}
\newcommand{\umaxb}
           {\bar u_{\textrm{\tiny max}}}
\newcommand{\alphag} 
           {{\alpha ^{\textrm{\tiny g}}}}
\newcommand{\alphagm} 
           {{\alpha ^{\textrm{\tiny g}}_{\textrm{\tiny max}}}}
\newcommand{\betag}  
           {\beta  ^{\textrm{\tiny g}}}
\newcommand{\rmin}{r^{\textrm{\rm\tiny min}}}
\newcommand{\rmax}{r^{\textrm{\rm\tiny max}}}
\newcommand{\rb}{r^{\textrm{\rm\tiny b}}}
\newcommand{\cmax}{c^{\textrm{\rm\tiny max}}}
\newcommand{\dmax}{d^{\textrm{\rm\tiny max}}}

\newcommand{\LLg} {\LL^{\textrm{\rm\tiny g}}}
\newcommand{\LLb} {\LL^{\textrm{\rm\tiny b}}}
\newcommand{\LLd} {\LL^{\textrm{\rm\tiny d}}}
\newcommand{\LLc} {\LL^{\textrm{\rm\tiny c}}}

\newcommand{\llg} {\ell^{\textrm{\rm\tiny g}}}
\newcommand{\llb} {\ell^{\textrm{\rm\tiny b}}}
\newcommand{\lld} {\ell^{\textrm{\rm\tiny d}}}
\newcommand{\llc} {\ell^{\textrm{\rm\tiny c}}}

\newcommand{\bll}  {\bar\ell}
\newcommand{\bllg} {\bar\ell^{\textrm{\rm\tiny g}}}
\newcommand{\bllb} {\bar\ell^{\textrm{\rm\tiny b}}}
\newcommand{\blld} {\bar\ell^{\textrm{\rm\tiny d}}}
\newcommand{\bllc} {\bar\ell^{\textrm{\rm\tiny c}}}

\newcommand{\tllg} {\tilde\ell^{\textrm{\rm\tiny g}}}
\newcommand{\tllb} {\tilde\ell^{\textrm{\rm\tiny b}}}
\newcommand{\tlld} {\tilde\ell^{\textrm{\rm\tiny d}}}
\newcommand{\tllc} {\tilde\ell^{\textrm{\rm\tiny c}}}

\newcommand{\MMg} {\MM^{\textrm{\rm\tiny g}}}
\newcommand{\MMb} {\MM^{\textrm{\rm\tiny b}}}
\newcommand{\MMd} {\MM^{\textrm{\rm\tiny d}}}
\newcommand{\MMc} {\MM^{\textrm{\rm\tiny c}}}

\newcommand{\Gg} {\Gamma^{\textrm{\rm\tiny g}}}
\newcommand{\Gb} {\Gamma^{\textrm{\rm\tiny b}}}
\newcommand{\Gd} {\Gamma ^{\textrm{\rm\tiny d}}}
\newcommand{\Gc} {\Gamma ^{\textrm{\rm\tiny c}}}

\newcommand{\Ng} {\mathbf{N}^{\textrm{\rm\tiny g}}}
\newcommand{\Nb} {\mathbf{N}^{\textrm{\rm\tiny b}}}
\newcommand{\Nd} {\mathbf{N}^{\textrm{\rm\tiny d}}}
\newcommand{\Nc} {\mathbf{N}^{\textrm{\rm\tiny c}}}
\newcommand{\tNg} {\mathbf{\tilde N}^{\textrm{\rm\tiny g}}}
\newcommand{\tNb} {\mathbf{\tilde N}^{\textrm{\rm\tiny b}}}
\newcommand{\tNd} {\mathbf{\tilde N}^{\textrm{\rm\tiny d}}}
\newcommand{\tNc} {\mathbf{\tilde N}^{\textrm{\rm\tiny c}}}

\renewcommand{\ng} {\mathbf{n}^{\textrm{\rm\tiny g}}}
\newcommand{\nb} {\mathbf{n}^{\textrm{\rm\tiny b}}}
\newcommand{\nd} {\mathbf{n}^{\textrm{\rm\tiny d}}}
\newcommand{\nc} {\mathbf{n}^{\textrm{\rm\tiny c}}}

\newcommand{\Mg} {{\mathbf{M}}^{\textrm{\rm\tiny g}}}
\newcommand{\Mb} {{\mathbf{M}}^{\textrm{\rm\tiny b}}}
\newcommand{\Md} {{\mathbf{M}}^{\textrm{\rm\tiny d}}}
\newcommand{\Mc} {{\mathbf{M}}^{\textrm{\rm\tiny c}}}

\newcommand{\mg} {{\mathbf{m}}^{\textrm{\rm\tiny g}}}
\newcommand{\mb} {{\mathbf{m}}^{\textrm{\rm\tiny b}}}
\newcommand{\md} {{\mathbf{m}}^{\textrm{\rm\tiny d}}}
\newcommand{\mc} {{\mathbf{m}}^{\textrm{\rm\tiny c}}}

\newcommand{\Ag} {A^{\textrm{\rm\tiny g}}}
\newcommand{\Ab} {A^{\textrm{\rm\tiny b}}}
\newcommand{\Ad} {A^{\textrm{\rm\tiny d}}}
\newcommand{\Ac} {A^{\textrm{\rm\tiny c}}}

\newcommand{\Atria} {A^{\vartriangle}}
\newcommand{\LLtria} {\LL^{\vartriangle}}
\newcommand{\lltria} {\ell^{\vartriangle}}
\newcommand{\blltria} {\bar\ell^{\vartriangle}}
\newcommand{\tlltria} {\tilde\ell^{\vartriangle}}
\newcommand{\Mtria} {{\mathbf{M}}^{\vartriangle}}
\newcommand{\mtria} {{\mathbf{m}}^{\vartriangle}}

\newcommand{\tMg} {\tilde M^{\textrm{\rm\tiny g}}}
\newcommand{\tMb} {\tilde M^{\textrm{\rm\tiny b}}}
\newcommand{\tMd} {\tilde M^{\textrm{\rm\tiny d}}}
\newcommand{\tMc} {\tilde M^{\textrm{\rm\tiny c}}}

\newcommand{\Zkf} {Z ^{k,f}}

\newcommand{\gammad} {\gamma^{\textrm{\rm\tiny d}}}
\newcommand{\gammab} {\gamma^{\textrm{\rm\tiny b}}}
\newcommand{\gammac} {\gamma^{\textrm{\rm\tiny c}}}

\renewcommand{\Lambda}  {\bar\lambda}
\newcommand{\Lambdad} {\bar\lambda^{\textrm{\rm\tiny d}}}
\newcommand{\Lambdab} {\bar\lambda^{\textrm{\rm\tiny b}}}
\newcommand{\Lambdac} {\bar\lambda^{\textrm{\rm\tiny c}}}
\renewcommand{\labelitemi} {\textbullet}
\renewcommand{\labelitemii}{$\diamond$}
\newcommand{\vois}[1]     {\stackrel{{#1}}{\sim}}
\definecolor{fyellow}{rgb}{1, 1, 0}


\title{A spatially explicit Markovian individual-based model for terrestrial plant dynamics\footnotetext{This work was partially supported by the French national research agency (ANR) within the SYSCOMM project ANR-08-SYSC-012 (MODECOL).}}
\author{Fabien Campillo\thanks{\protect\url{Fabien.Campillo@inria.fr} --- 
           Project--Team MERE, INRIA/INRA, UMR ASB, bât. 29, 2 place Viala, 34060 Montpellier cedex 06, France} 
   \and Marc Joannides\thanks{\protect\url{marc.joannides@univ-montp2.fr} ---
   Université Montpellier 2 / I3M, case courrier 51,
        place Eugène Bataillon, 34095 Montpellier cedex 5; this author is associate researcher for Project--Team MERE, INRIA/INRA, UMR ASB.}
      }
\date{April 23, 2009}

\maketitle

\begin{abstract}
An individual-based model (IBM) of a spatiotemporal terrestrial ecological population is proposed. This model is spatially explicit and features the position of each individual together with another characteristic, such as the size of the individual, which evolves according to a given stochastic model. The population is locally regulated through an explicit competition kernel.  The IBM is represented as a measure-valued branching/diffusing stochastic process. The approach allows \fenumi\ to describe the associated Monte Carlo simulation and \fenumii\ to analyze the limit process under large initial population size asymptotic. The limit macroscopic model is a deterministic integro-differential equation.



\end{abstract}



\section{Introduction}

Our aim is to present a spatially explicit individual-based model of a spatiotemporal terrestrial ecological system. We consider a family of individuals whose state includes their position and another characteristic such as their size. During the life of an individual, its position remains constant and its characteristic evolves according to a given stochastic continuous model. Each individual is subject to random punctual mechanisms: natural mortality, mortality due to competition, and reproduction. The individuals interact through the competition mechanism.

So-called individual-based models (IBM) are not new in the field of theoretical ecology where they have become more widespread thanks to the development of computers. Since the publication of the significant paper by Huston {\it et al} \cite{huston1988a}, many studies have been devoted to these models \cite{deangelis1992a}, \cite{grimm2005a}. IBMs appear mainly as computer simulators based on empirical rules. IBMs are modeling tools that require design work, and which cannot be reduced to the indiscriminate use of empirical rules. The study of the behavior and properties of IBMs still requires mathematical analyses.

Mathematical representations of such models in continuous time and space was introduced by Bolker-Pacala \cite{bolker1997a} and Dieckmann-Law \cite{dieckmann2000b}. These authors also derived a truncated moment method that determines the time evolution of the approximated first moments. Méléard and co-workers \cite{fournier2004a}, \cite{champagnat2006a}, \cite{champagnat2006b}, \cite{champagnat2007b} proposed a rigorous mathematical setup for these models: the dynamics of the ecological system are described as a measure-valued branching/diffusing random process. The limit law behavior of this type of microscopic Markov process, as certain parameters such as the size of the initial population tend to infinity, is relatively well known. The macroscopic limit models can be of different nature: deterministic, such as ordinary, integro, or partial differential equations; or probabilistic, such as stochastic partial differential equations or superprocesses  \cite{ethier1986a}, \cite{etheridge2000a}. 

Within this framework, various kinds of dynamics have been considered: terrestrial plants \cite{fournier2004a}, Darwinian evolution \cite{champagnat2006a,champagnat2006b, champagnat2007b}, phytoplankton aggregation \cite{elsaadi2006b}, age-structured population \cite{meleard2008a}. 

Méléard and co-workers also described an exact Monte Carlo procedure to simulate the associated microscopic stochastic process. These algorithmic aspects which are very important in practice, deserve specific attention.

\medskip

The objective is threefold. Firstly, we set out an IBM mathematical framework for ecosystems such as terrestrial plant systems. Secondly, we develop the associated Monte Carlo algorithm. Thirdly, we derive the macroscopic behavior of the IBM. 

We focus on terrestrial plant ecosystem dynamics models that are spatially explicit with an explicit representation of the competition interactions between individuals. This is one of the most active areas of computational ecology \cite{gratzer2004a}, \cite{neuhauser2001a}, \cite{berger2008a}.

In \cite{fournier2004a} the authors present a model for terrestrial plant ecosystem dynamics where they consider only the position of the plant individuals. In our work we extend this model further to include a continuously evolving characteristic such as the size of the individuals.

We describe the Bolker-Pacala-Dieckman-Law model in Section 2. The Monte Carlo simulation method is presented in Section 3. An example of a terrestrial plant ecosystem is presented in Section 4. The Markov representation of the IBM is described in Section 5. Finally the large population limit is analyzed in Section 6. The generic mathematical approach used to prove the weak convergence of measure-valued process in Section 6 is detailed in the Appendix.

\section{The model}

We consider a family of individuals that live in a set of the form:
\[
    \XX = \DD \times \R^d
\]
where $\DD$ is a measurable connected subset of $\R^n$. The state $x=(p,r)\in\XX$ of an individual represents its position $p$ in the physical space $\DD$ and an associated characteristic $r$ that could represent its size or its maturation age.

It is convenient to represent an individual at point $x\in\XX$ as the Dirac measure $\delta_x$; hence the population at time $t$ will be:
\begin{equation}
\label{eq.nu}
  \nu_t(\rmd x) = \sum_{i=1}^{N_t} \delta_{x^i_t}(\rmd x)
\end{equation}
where the sum is over all individuals alive at time $t$ and $N_t$ is the size of the population at time $t$.

Let $\MM_F(\XX)$ be the set of finite positive measures on $\XX$, and $\MM$ be the set of finite point measure on $\XX$, i.e.
\[
  \textstyle
  \MM \eqdef \big\{\sum_{i=1}^N\delta_{x^i}\,;\,N\in\N^*,\, x_i\in\XX \big\}\,.
\]
The duality operator between the measures and the functions will be denoted:
\[
  \crochet{\mu,f}=\int_\XX f(x)\,\mu(\rmd x)
\]
(note that $\bcrochet{Z}_t$ will denote the quadratic variation of a process $Y_t$). 
Hence, $(\nu_t)_{t\geq 0}$ is a $\MM$-valued process and, according to (\ref{eq.nu}):
\[
  N_t \eqdef \crochet{\nu_t,1}
\]
is the size of the population at time $t$. By abuse of notation, ``$x\in\nu_t$'' will specify that an individual in the state $x$ belongs to the family $\nu_t$ at time $t$. 

\begin{remark}[numbering convention]
\label{rem.conv.numerotation}
In practice, starting from a family labelled from $1$ to $N$, we use the following numbering of particles: \fenumi\ in case of birth, the new individual will be labelled $N+1$; \fenumii\ in case of death of the individual $i$, the first $(i-1)^\textrm{th}$ labels remains unchanged and the last $(N-i)^\textrm{th}$ labels are left shifted (i.e. $j\to j-1$). As pointed out in \cite{fournier2004a} and \cite{champagnat2006b}, this numbering convention has no influence on the law of the process we will describe, it affects only the trajectorial realizations of the process.
\end{remark}

\bigskip
	
Considering the state $\nu = \sum_{i=1}^N\delta_{x^i}$ of the family
at a given time, an individual in state  $x\in \nu$ will be subject to
3 types of \emph{punctual} events occurring at specific rates:
\begin{quote}
\begin{description}

\item[Intrinsic death:] This individual disappears at a rate $\lambdad(x)$ which 
   may depend on its state. This death is called ``intrinsic''
   as it does not depend on the state of all the population $\nu$. It
   represents the ``natural death'' as opposed to ``competition death''.

\item[Competition death:] This individual disappears at a rate $\lambdac(x,\nu)$
  which may depend on its state $x$ and on the state $\nu$ of the population. 
  We suppose that $\lambdac(x,\nu)$ is of the form:
  \begin{align}
  \label{eq.comp.kernel}
    \lambdac(x,\nu) 
    = \sum_{y\in\nu} u(x,y) 
    = \int_{\XX} u(x,y)\,\nu(\rmd y)\,.
  \end{align}
  The competition kernel $u(x,y)$ is the contribution of an individual located at $y$ to the
  competition affecting an individual located at $x$.   
\item[Birth and dispersal:] This individual gives birth to a new individual 
  at a rate $\lambdab(x)$ which may depend on its state. The state $y\in\XX$ of the new 
  individual will be determined by a given dispersal kernel  (see Remark \ref{rem.kernel} later).

\end{description}
\end{quote}
Between discrete events of birth or death, the size of the population remains unchanged as well as the position of the individuals, and the population state $(x^i_t)_{1\leq i\leq N}$ is subject to a \emph{continuous} mechanisms:
\begin{quote}
\begin{description}
\item[Displacement:]  Over time, the characteristic component of each individual $i$ evolves 
   in the state space $\R^d$ in interaction with the evolution of all other
   individuals according to the following system of stochastic differential
   equations (SDE):
 	\begin{align}
	\label{eq.growth.tilde}
  		\rmd
  		\left(\begin{smallmatrix}p^i_t \\ r^i_t \end{smallmatrix}\right)
  		=
  		\left(\begin{smallmatrix}0 \\ \tilde g(x^i_t,\nu_t) \end{smallmatrix}\right)
 		 \,\rmd t
 		 +
 		 \left(\begin{smallmatrix}0  \\ \tilde\sigma(x^i_t,\nu_t) \end{smallmatrix}\right)
 		 \,\rmd \mathbf{B}^i_t
	\end{align}
	where $(\mathbf{B}^i_t)_{t\geq 0}$ are independent standard Brownian motions.
	To simplify the notation, Equation (\ref{eq.growth.tilde}) will be represented as: 
	\begin{align}
	\label{eq.growth}
	  \rmd x^i_t = g(x^i_t,\nu_t)\,\rmd t
	               +\sigma(x^i_t,\nu_t)\,\rmd \mathbf{B}^i_t\,.
	\end{align}
Let
\begin{align*}
   a(x,\nu) \eqdef  \sigma(x,\nu) \, \sigma^*(x,\nu)\,.
\end{align*}
We define the
associated flow operator:
\begin{align*}
   \nu_t
   =
   \bbF(t,s;\nu_s)
\end{align*}
defined for all $s\leq t$ between two successive punctual events (i.e. between two successive jump of the population size). Note that the initial condition distribution law in SDE (\ref{eq.growth.tilde}) or (\ref{eq.growth}) is handled by the the dispersal kernel (see Remark \ref{rem.kernel} later).
\end{description}

\end{quote}
We suppose that these four mechanisms and the dispersal mechanism are mutually independent.

\begin{remark}[dispersal kernel]
\label{rem.kernel}
An individual in state  $x=(p_x,r_x)$ will give birth to a new individual in state
$y=(p_y,r_y)=x+y'=(p_x+p_{y'},r_x+r_{y'})$ (``$y$'' and ``${y'}$'' will denote respectively the absolute state and the relative state of the new individual). The state of the new individual is given by a relative kernel.

In the present application it will be convenient to consider a ``mixed relative/absolute'' formulation: an individual in state  $x=(p_x,r_x)$ will give birth to a new individual in state
$(p_x+p_z,r_z)$ according to a kernel $D(x,\rmd z)$ with $z=(p_z,r_z)$. We suppose that this kernel admits a density:
\begin{align}
\label{eq.D.mixed}
  D(x,\rmd z)
  =
  D(x,z)\,\rmd z\,.
\end{align}
This formulation is natural: the position $p_x+p_z$ of the new individual will be relative to $p_x$ and its characteristic $r_z$ will be absolute. For phenotypic trait dynamics \cite{champagnat2007b}, it is more natural to consider a relative mutation kernel for the $r$-component. 
\end{remark}

\section{Monte Carlo simulation}
\label{sec.monte.carlo}

We now describe the dynamic of the process starting from an initial population state $\nu$.
Independently of each other, an individual with state $x$ in the population $\nu$ has three independent exponential clocks that control the occurrence of the events: \fenumi\ a birth clock with rate $\lambdab(x)$, \fenumii\ an intrinsic death clock with rate $\lambdad(x)$,
\fenumiii\ a competition death clock with rate $\lambdac(x,\nu)$.
\begin{enumerate}
\item
When the birth clock rings, the individual $x$ gives birth to a new individual with a state
$z\in\XX$ determined by the dispersal kernel $D(x,\rmd z)$:
\[
   \nu \to \nu + \delta_{(p_x+p_z,r_z)}\,.
\]
\item
When the intrinsic or competition death clock rings, the individual $x$ is removed 
from the population.
\[
   \nu \to \nu - \delta_x\,.
\]

\end{enumerate}
Between any birth or death event, the state of all the population evolves according to  (\ref{eq.growth}) which corresponds to a system of $N$ interacting SDE's ($N=\crochet{\nu,1}$ is the size of the population).

\bigskip

Considering individual clocks is cumbersome, a more efficient Monte Carlo procedure will rely on the existence of a \emph{global clock} that dominates all punctual phenomena (birth, natural death, competition death). That existence holds true when all the different local clocks are uniformly bounded: then, given an individual chosen at random in the population, the type of punctual phenomenon to be considered is determined by a sampling technique, and it is decided whether the chosen phenomenon is actually applied or not by an acceptance/rejection sampling technique. The existence of a uniform bound avoids explosion phenomena, i.e. accumulation of infinitely many events at a given time.
\begin{hypotheses}
\label{hyp.horloge.globale}
We suppose that there exists positive real numbers $\lambdabm$,  $\lambdadm$ and $\umax$ such that:
\begin{align*}
   \lambdab(x) &\leq \lambdabm\,,
   &
   \lambdad(x) &\leq \lambdadm\,,
   &
   u(x,y) &\leq \umax\,.
\end{align*}
Hence:
\[
   \lambdac(x,\nu)
   \leq 
   \umax\,\crochet{\nu,1}\,.
\]
We also suppose that the dispersal kernel (\ref{eq.D.mixed}) satisfies:
\begin{align}
\label{eq.majoration.D}
   D(x,z) &\leq \kappa\,\bar D(z)
\end{align}
where $\bar D$ is a probability density function.
\vskip-1.3em\carre
\end{hypotheses}

\medskip

Let $T_0=0$, and suppose $T_{k-1}$ and $\nu_{T_{k-1}}$ given. We describe now how to simulate $\nu_{T_{k}}$ starting from $\nu_{T_{k-1}}$. In order to determine the instant $T_{k}$ where the next event \emph{could} take place, we should bound the different rates uniformly in space but also in time. This is possible thanks to the Hypothesis \ref{hyp.horloge.globale}. From the instant $T_{k-1}$ to the instant $T_{k}$ of the next event, i.e. along the time interval $[T_{k-1},T_{k}[$, the population size is $N = \crochet{\nu_{T_{k-1}},1}$.

At the scale of the population, the maximum rate of events (birth, natural death, death by competition) is bounded by:
\begin{align}
\label{eq.taux.max.v1}
  \gamma
  \eqdef
  \gammab+\gammad+\gammac
  \quad
  \textrm{with }
  \left\{\begin{array}{l}
    \gammab \eqdef \kappa\,\lambdabm \,N \\
    \gammad \eqdef\lambdadm \, N \\
    \gammac \eqdef\umax \,N^2
  \end{array}\right.
\end{align}
where
\begin{quote}
\begin{description}
\item
[$\gammab$] is an upper bound of the birth rate at the population scale (i.e. the birth rate at the population scale if all individuals have the same birth rate $\lambdabm$).
\item
[$\gammad$] is an upper bound of the natural death rate at the population scale (i.e. the natural death rate at the population scale if all individuals have the same natural death rate $\lambdadm$).
\item
[$\gammac$] is an upper bound of the rate of death by competition at the population scale.
\end{description}
\end{quote}
An acceptance/rejection method will permit us to correct the fact that these three terms are upper bounds for the actual rates.

\paragraph{Iteration $\nu_{T_{k-1}} \to \nu_{T_{k}} $:}  
\label{algo.simu.1.detail}
\begin{quote}
\it 
\begin{enumerate}

\item Let $N=\crochet{\nu_{T_{k-1}},1}$ be the population size.

\item
Computation of the global rate $
  \gamma
  \eqdef
  \gammad+\gammab+\gammac
$
with (\ref{eq.taux.max.v1}).

\item
Simulation of the next event instant:
\[
   T_k = T_{k-1}+S_k\,,
   \quad\textrm{with }S_k\sim\textrm{Exp}(\gamma)\,.
\]

\item
Computation of the system evolution  between the two instants:
\[
   \nu_{T_k^-}
   =
   \bbF(T_k,T_{k-1},\nu_{T_{k-1}})\,.
\]
In practice, the system is simulated with an  Euler discretization scheme.

\item
Chose $x$ at random uniformly in $\nu_{T_k^-}$; chose at random the nature of the next event according the probability values $(\gammab/\gamma, \gammad/\gamma, \gammac/\gamma)$:
\begin{quote}
\begin{itemize}
\item birth: choose $z'$ according to the law $\bar D(z)\rmd z$ and let
\[
   \nu_{T_k} 
   = 
   \left\{\begin{array}{ll}
     \nu_{T_{k}^-}+\delta_{(p_x+p_{z'},r_{z'})} 
       & \textrm{with probability } \frac{\lambdab(x)\,D(x,z')}{\lambdabm \,\kappa\,\bar D(z')}
    \\
     \nu_{T_{k}^-}
       & \textrm{with probability } 1-\frac{\lambdab(x)\,D(x,z')}{\lambdabm \,\kappa\,\bar D(z')}
   \end{array}\right.
\]

\item natural death: 
\[
   \nu_{T_k} 
   = 
   \left\{\begin{array}{ll}
     \nu_{T_{k}^-} - \delta_x 
       & \textrm{with probability } \frac{\lambdad(x)}{\lambdadm }
    \\
     \nu_{T_{k}^-}
       & \textrm{with probability } 1-\frac{\lambdad(x)}{\lambdadm }
   \end{array}\right.
\]

\item competition death:
chose $y$ at random uniformly in $\nu_{T_k^-}$ and let
\[
   \nu_{T_k} 
   = 
   \left\{\begin{array}{ll}
     \nu_{T_{k}^-}- \delta_x 
       & \textrm{with probability } \frac{u(x,y)}{\umax}
    \\
     \nu_{T_{k}^-}
       & \textrm{with probability } 1-\frac{u(x,y)}{\umax}
   \end{array}\right.
\]
\end{itemize}
\end{quote}
\vskip-1em
\carre

\end{enumerate}
\end{quote}
The algorithm is detailed in Algorithm~\ref{algo.simu.1}.

We will present numerical tests in a future companion article \cite{campillo2009d} where we will also propose other versions of this algorithm.

\begin{algorithm}[p]
\caption{\sl Simulation algorithm, first version (see description in Section \ref{sec.monte.carlo}).}
\label{algo.simu.1}
\begin{center}
\begin{minipage}{12cm}
\hrulefill\\[-1em]
\mbox{}
\algsetup{linenosize=\tiny}
\algsetup{linenodelimiter=}
\begin{algorithmic}[1]
\STATE $T_0\ot 0$, $k\ot 0$, initialization of $\nu_{T_0}$
\WHILE {$T_k\leq T_{\textrm{\tiny max}}$}
	\STATE  $k\ot k+1$
    \STATE  $N\ot \crochet{\nu_{T_{k-1}},1}$
    \STATE  $\gammab \ot \kappa\,\lambdabm \, N$,  
            $\gammad \ot \lambdadm \, N$,
            $\gammac \ot \umax \, N^2$,
            $\gamma \ot \gammab + \gammad + \gammac$
    \STATE  $S\sim \text{Exp}(\gamma)$ , $T_k\ot T_{k-1}+S$
    \STATE  $\nu_{T_k^-} \ot    \bbF(T_k,T_{k-1},\nu_{T_{k-1}})$
	\COMMENT{evolution (Euler scheme)}

    \STATE  $u\sim\mathcal{U}[0,1]$
    \STATE  choose $x$ at random uniformly in $\nu_{T_k^-}$

    \IF {$u\in [0,\gammab/\gamma]$}
    		\STATE $v\sim\mathcal{U}[0,1]$
      	\STATE $z'\sim \bar D(z)\,\rmd z$
    		\IF { $v \leq (\lambdab(x)\,D(x,z'))/(\lambdabm\,\kappa\,\bar D(z'))$}
      		\STATE $\nu_{T_k}
		         \ot\nu_{T_{k^-}} +\delta_{(p_x+p_{z'},r_{z'})}$
			\COMMENT{birth}
    		\ENDIF 

    \ELSIF{$u\in [\gammab/\gamma,(\gammab+\gammad)/\gamma]$}
    		\STATE $v\sim\mathcal{U}[0,1]$
    		\IF { $v \leq \lambdad(x)/\lambdadm$}
      		\STATE $\nu_{T_k}\ot\nu_{T_{k^-}} - \delta_x$
			\COMMENT{natural death}
    		\ENDIF 

    \ELSE
	    \STATE choose $y$ at random uniformly in $\nu_{T_k^-}$
	    \STATE $v\sim \mathcal{U}[0,1]$
	    \IF {$v \leq u(x,y)/\umax$}
	    		\STATE $\nu_{T_k}
			     \ot\nu_{T_{k^-}}- \delta_x$
			\COMMENT{competition death}
		\ENDIF
  	\ENDIF

\ENDWHILE
\end{algorithmic}
\hrulefill
\end{minipage}
\end{center}
\end{algorithm}

\section{A forest dynamic model}
\label{sec.foret}

\subsection*{Competition model: zone of influence approach}

\begin{figure}
\begin{center}
\setlength{\fboxsep}{0pt}
\includegraphics[width=7cm]{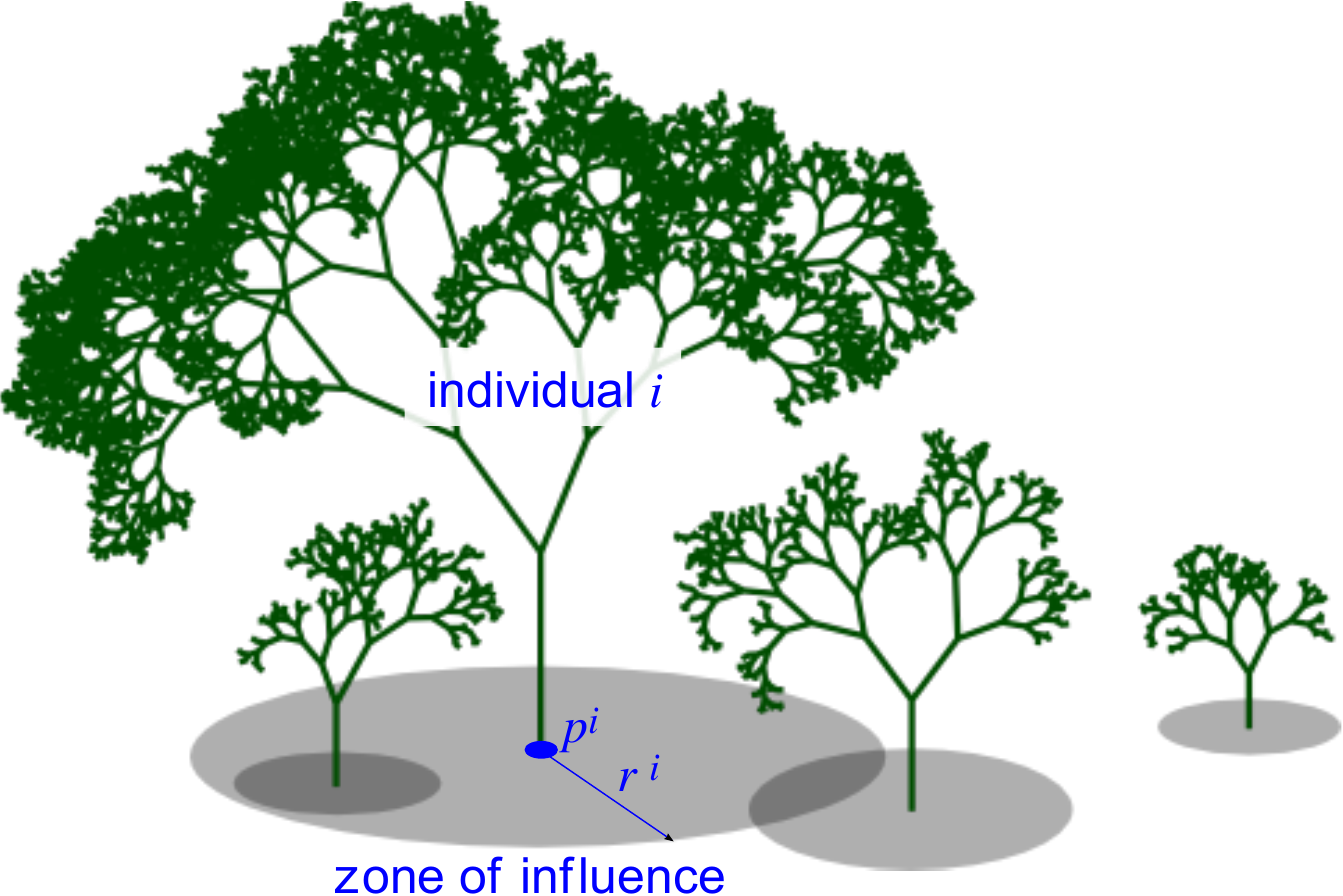}\hskip2cm
\includegraphics[width=4cm]{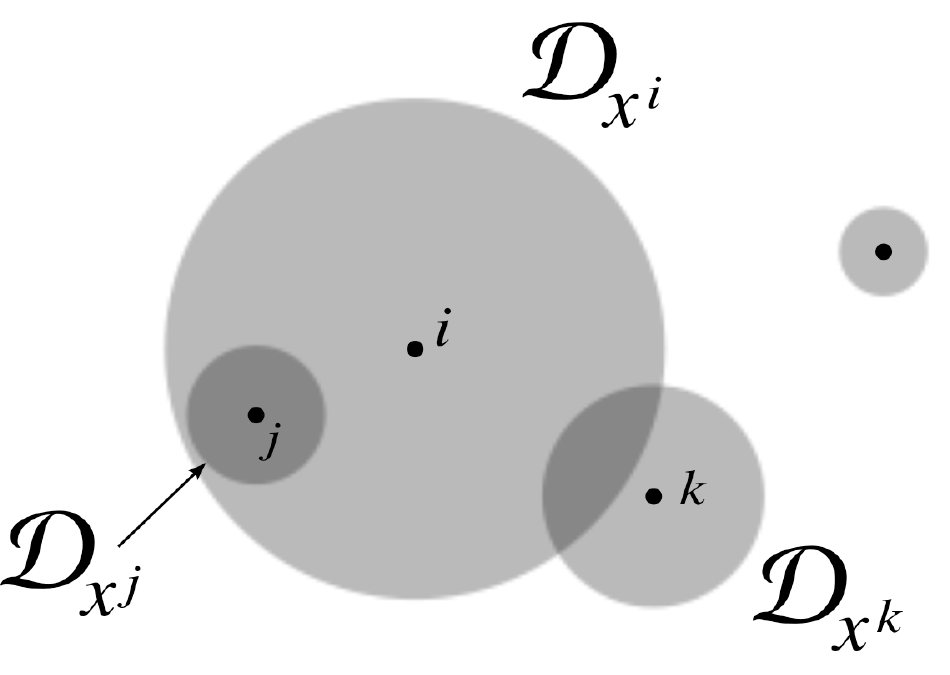}
\end{center}
\caption{\it 
{\bfseries Left:} Each individual tree $i$ is associated with a state $x^i=(p^i,r^i)$
   where $p^i$ is its position in $\DD\subset \R^2$ and $r^i\in
  [\rmin,\rmax]\subset \R_+$ is the radius of its zone of influence. This area represents 
  the portion of ground which the individual needs to ensure its growth.
  {\bfseries Right:} The more the zone of influence $\DD_x$ of an individual in state $x$ intersects with other individuals ZOI, the higher the strength of competition is, and the greater the associated death rate is. The ZOI appears in the definition of the competition kernel $u(x,y)$ in (\ref{eq.lambda.c}) and in the growth model (\ref{modele.richards}).} 
\label{fig.trees}
\end{figure}

Here we consider a population of trees. The state of each individual tree $i$ is:
\[
   x^i 
   \eqdef 
   (p^i,r^i)
   \in 
   [0,L]^2\times [\rmin,\rmax]\subset \R^2\times \R_+
\]
where $p^i$ is the position of the tree in a parcel $\DD=[0,L]^2$ and $r^i$ is the radius of its \emph{zone of influence} (cf. Figure \ref{fig.trees}). This zone of influence is the disk centered in $p^i$ with radius $r^i$ which symbolizes the portion of the ground that the individual needs to ensure its growth. For $x=(p,r)\in\nu$ let: 
\[
  \DD_x \eqdef \textrm{disk of center $p$ and radius $r$.}
\]

The local interaction kernel $u(x,y)$ introduced in (\ref{eq.comp.kernel}) is of the form:
\begin{align}
\label{eq.lambda.c}
  u(x,y)
  =
  \left\{\begin{array}{ll}
    \umax\,\frac{\textrm{Area}(\DD_x\cap \DD_y)}{\textrm{Area}(\DD_x)}
    &\textrm{if }x\neq y\,, 
    \\
    0 & \textrm{otherwise.}
  \end{array}\right.
\end{align}
The surface area $\textrm{Area}(\DD_x)$ of the zone of influence associated with an individual in state $x$ in $\nu$ represents the amount of resources (e.g. sunlight, water, nutrients...) needed for growth; $u(x,y)$ is the strength of competition experienced by the individual in state $x$ from an individual in state $y$ in the population $\nu$.

The computation of the areas in (\ref{eq.lambda.c}) could be cumbersome for large population sizes. 
However, there are many alternatives \cite{berger2008a}.

\subsection*{Birth and dispersal model}

Birth occurs at a rate $\lambdab(x)$ defined by:
\begin{align*}
   \lambdab(x)
   =
   \lambdab(p,r)
   =
   \lambdabm\,\frac{r}{\rmax}\,\indic_{\{r\geq \rb\}}
   \leq
   \lambdabm\,.
\end{align*}
This birth rate can be understood as a fertility model: only individual with ZOI radius
greater than $\rb$ can give birth; and the greater this radius is, the greater the rate is.

For the dispersal kernel $D(x,z)$, we can consider two possibilities:
\begin{enumerate}

\item
{\it ``Parcel in forest model''} -- Here $\DD$ is the torus $[0,L]^2$, and we consider an homogeneous kernel:
\begin{align*}
  D(x,z)
  =
  D(z)
  =
  D_1(p_z)\, D_2(r_z)\,.
\end{align*}
In this example we consider a slightly different case where the initial condition on $r_z$ does not admit a density:
For example:
\begin{align*}
   D(x,\rmd z)=D(x,\rmd p_z\times\rmd r_z)
   &=
   \bar\NN(0,\sigma^2\,I;\rmd p_z) 
   \times
   \delta_{\rmin}(\rmd r_z) 
\end{align*}
where $\bar\NN$ is a ``Gaussian law on the torus $\DD$'' (\footnote{That is $\bar\NN(0,\sigma^2\,I;B) \eqdef \P(X\in \cup_{p,q\in\Z} B+(p\,L,q\,L))$ for any Borelian set $B$ of $[0,L[^2$ where  $X\sim \NN(0,\sigma^2\,I)$.
}).

The condition (\ref{eq.majoration.D}) is fulfilled but in fact the present setup is even simpler. Indeed, the part {\fenumv} of the algorithm proposed at the end of section \ref{sec.monte.carlo} is now:

\begin{quote}
\begin{itemize}
\item birth: 
\[
   \nu_{T_k} 
   = 
   \left\{\begin{array}{ll}
     \nu_{T_{k}^-}+\delta_{(p_x+p_{z'},r_{z'})} 
       & \textrm{with probability } \frac{\lambdab(x)}{\lambdabm}
    \\
     \nu_{T_{k}^-}
       & \textrm{with probability } 1-\frac{\lambdab(x)}{\lambdabm}
   \end{array}\right.
\]
where $z$ is simulated according to the law $D_1(p_z)\, D_2(r_z)$.
\end{itemize}
\end{quote}
This setup is periodic and could illustrate the case of a squared parcel imbedded in a forest: the descendants of the individual of the parcel may integrate the contiguous parcels and the individuals in the contiguous parcels may integrate the parcel under consideration.

\item 
{\it``Island model''} --
Here $\DD=[0,L]^2$ is a bounded squared parcel and we suppose that when an individual is closed to the border of $\DD$ a portion of its offspring is lost in the water and do not integrate the population. Hence, the individuals near the border have a lower effective fertility rate compared to the individuals closed to the center of the parcel. 

For simplicity we suppose that the birth fertility rate is constant for individuals at distance greater than $r_0$ from the border, for these individuals suppose that the dispersal kernel is homogeneous: $D_1(p_x,p_z)=\bar D_1(p_z)$ (with the support of $\bar D_1(p_z)$ is included in a disk of center 0 and radius $r_0$), suppose that all these individuals have a constant fertility rate $\lambdab$. For individual located in $p_x$ at most $r_0$ from the border, then 
\[
  D_1(p_x,p_z)=C_x\times\indic_{[0,L]^2}(p_x+p_z) \times \bar D_1(p_z)
\]
with $C_x = [\int_{[0,L]^2} \bar D_1(p_y-p_x)\,\rmd p_y]^{-1}$. Here we suppose that the fertility rate associated with the later points is non-constant and equal to $\lambdab(x) = \lambdab/C_x$, i.e. the more the kernel $D_1$ intersects the border, the lower the birth rate is.

\end{enumerate}

In the previous model we compensated the birth rate for individual close to the border. Without this mechanism we get a model where there is an accumulation of new individuals on the border of the parcel.  In a way, the ``parcel in forest'' model is the more realistic and therefore, this is the one we used for simulation purposes.

\subsection*{Growth model}

Suppose that the radius of the ZOI is solution of the deterministic equation:
\begin{align}
\label{eq.ri}
	  \dot r^i_t = g(x^i_t,\nu_t)
\end{align}
while the position remains unchanged. This last equation is coupled with 
$\dot p^i_t = 0$ so that we get a system of ODE's for $x^i_t$. For (\ref{eq.ri}), we consider a model of the form:
\begin{subequations}
\label{modele.richards}
\begin{align}
\label{modele.richards.0}
  \dot r^i_t 
  &=
  \psi(x^i_t,\nu_t)\,  R(r^i_t) 
  \,,\quad
  r_0=\rmin
\end{align}
where $R(r^i_t)$ corresponds to a standard growth model such as the Richards model \cite{damgaard2004b}: 
\begin{align}
\label{modele.richards.1}
  R(r^i_t) 
  &\eqdef
  \alphagm \, \frac{1}{1-\betag}\,r^i_t\,
  \Big[\big(\frac{r^i_t}{\rmax}\big)^{\betag-1} - 1\Big]
\end{align}
with $\betag\neq 1$ and $\psi(x^i_t,\nu_t)$ is an expression between 0 and 1. The case $\psi(x^i_t,\nu_t)=1$ corresponds to the best condition for growth: in this case its growth is described by the Richards model $\dot r^i_t = R(r^i_t)$. The smaller $\psi(x^i_t,\nu_t)$ is, the more the growth conditions of the individual $i$ are degraded. Small $\psi(x^i_t,\nu_t)$ corresponds to the situation where the individual $i$ is surrounded by many other individuals. We may think of:
\begin{align}
\label{modele.richards.2}
  \psi(x^i_t,\nu_t)
  &=
  \Big[1-C_{\textrm{\rm g}} \, \lambdac(x^i_t,\nu_t)\Big]^+
  \in[0,1]
\end{align}
or, more generally, $\psi(x^i_t,\nu_t) = 
\Psi(\lambdac(x^i_t,\nu_t))$ where
$\Psi:\R_+ \mapsto [0,1]$ is any continuous, decreasing function such that $\Psi(0)=1$.
 
\end{subequations}



\section{Markov representation of the process $(\nu_t)_{t\geq 0}$}

\subsection{Identification of the infinitesimal generator}
\label{sub.sec.infinitesimal.generator}

We introduce the following set $\DD$ of test functions $\Phi:\MM_F(\XX)\mapsto \R$ of the form:
\[
   \Phi(\nu) = F(\crochet{\nu,f})
\]
for any function $f:\XX\mapsto\R$ and $F:\R\mapsto\R$ twice continuously differentiable, bounded with bounded derivatives.

At the end of this we will present a particular case of such test functions $\Phi$.

\begin{lemma}[expression for $\Phi(\nu_t)$]
\label{lemma.Phi.nu}
For any $\Phi=(F,f)\in\DD$, $\Phi(\nu_t)=F(\crochet{\nu_t,f})$ satisfies:
\begin{align*}
  \Phi(\nu_t)
  &=
  \Phi(\nu_0)
  +
  \int_0^t \LL \Phi(\nu_{s})\;\rmd s+ \mathbf{M}_{\Phi,t}(\nu)
\\
  &=
  \Phi(\nu_0)
  +
  \sum_{\vartriangle=\textrm{\rm d,b,c,g}}
    \Big\{
       \int_0^t \LLtria \Phi(\nu_{s})\;\rmd s
       +
       \Mtria_{\Phi,t}(\nu)
    \Big\}
\end{align*}
where
\begin{enumerate}

\item
$\LL=\LLd+\LLb+\LLc+\LLg$ are the infinitesimal generator defined by:
\begin{subequations}
\label{eq.L}
\begin{align}
\label{eq.Ld}
  \LLd \Phi(\nu)
  & \eqdef
  \int_{\XX} \lambdad(x)\,
     [\Phi(\nu-\delta_{x}) - \Phi(\nu)]\,
    \nu(\rmd x)\,.
\\
\label{eq.Lb}
  \LLb \Phi(\nu)
  &\eqdef
  \int_\XX
       \lambdab(x)\,
       \Big[ \int_{\XX}  
                [\Phi(\nu+\delta_{(p_x+p_z,r_z)})-\Phi(\nu)]
       \,D(x,z)\,\rmd z\Big]
   \;\nu(\rmd x)
\\
\label{eq.Lc}
  \LLc\Phi(\nu)
  &\eqdef
  \int_\XX 
     \Big[\int_\XX
     u(x,y)\,\,\nu(\rmd y)
     \Big]
     [\Phi(\nu-\delta_{x}) - \Phi(\nu)]\,
     \;\nu(\rmd x)
\\
\label{eq.Lg}
  \LLg\phi(\nu)
  & \eqdef
  \vphantom{\int_\XX}
  F'(\crochet{\nu,f})
  \times\textstyle
  \crochet{
          \nu
          ,
          \GG f(\cdot,\nu)
  }
  +
  \demi
  F''(\crochet{\nu,f})
  \times
  \crochet{
      \nu
           ,
           |\nabla f\cdot \sigma(\cdot,\nu)|^2
  }
\end{align}
\end{subequations}
corresponding respectively to natural death, birth, death by competition and growth respectively, and with:
\begin{align*}
  \GG f(x,\nu)
  \eqdef\textstyle
  \nabla f(x) \cdot g(x,\nu)
         + \demi\,\frac{\partial^2 f(x)}{\partial x_\ell\,\partial x_{\ell'}}
             \,a_{\ell\ell'}(x,\nu)
\end{align*}

\item
the martingale term $\mathbf{M}_{\Phi,t}(\nu)=\Md_{\Phi,t}(\nu) + \Mb_{\Phi,t}(\nu) + \Mc_{\Phi,t}(\nu) + \Mg_{\Phi,t}(\nu)$ is the sum of four martingales defined by:
\begin{subequations}
\label{eq.M}
\begin{align}
\nonumber
  \Md_{\Phi,t}(\nu)
  & \eqdef
  \int_0^t\int_{\N^*}\int_0^1 
     \indic_{(i\leq N_{s^-})}\,
     \indic_{(\theta \leq \lambdad(x^i_{s^-})/\lambdadm)}\,
\\
\label{eq.Md}
  &\qquad\qquad \qquad\qquad
    \times
       [\Phi(\nu_{s^-}-\delta_{x^i_{s^-}}) - \Phi(\nu_{s^-})]\,
     \tNd(\rmd s,\rmd i,\rmd \theta)\,,
\\
\nonumber
  \Mb_{\Phi,t}(\nu)
  & \eqdef
  \int_0^t\int_{\N^*}\int_{\XX}\int_0^1 
     \indic_{(i\leq N_{s^-})}\,
     \indic_{\{\theta \leq (\lambdab(x^i_{s^-})\, D(x^i_{s^-},z))/(\lambdabm\,\kappa\,\bar D(z))\}}\,
\\
\label{eq.Mb}
  &\qquad\qquad \qquad\qquad
    \times
     [\Phi(\nu_{s^-}+\delta_{(p_{x^i_{s^-}}+p_z,r_z)}) - \Phi(\nu_{s^-})]\,
     \tNb(\rmd s,\rmd i,\rmd z,\rmd \theta)\,,
\\
\nonumber
  \Mc_{\Phi,t}(\nu)
  & \eqdef
  \int_0^t\int_{\N^*}\int_{\N^*}\int_0^1 
     \indic_{(i\leq N_{s^-})}\,
     \indic_{(j\leq N_{s^-})}\,
     \indic_{(\theta \leq u(x^i_{s^-},x^j_{s^-})/\umax)}\,
\\
\label{eq.Mc}
  &\qquad \qquad \qquad\qquad
    \times
     [\Phi(\nu_{s^-}-\delta_{x^i_{s^-}}) - \Phi(\nu_{s^-})]\,
     \tNc(\rmd s,\rmd i,\rmd j,\rmd \theta)
     \,,
\\
\label{eq.Mg}
  \Mg_{\Phi,t}(\nu)
  &\eqdef
  \int_0^t 
        \textstyle
        F'(\crochet{\nu_s,f})\times
        \sum_{i=1}^{N_s}\nabla f(x^i_s)\cdot \sigma(x^i_s,\nu_s)\,\rmd \mathbf{B}^i_s
        \,,
\end{align}
\end{subequations}
and
\begin{enumerate}

\item
$\Nd$ is a Poisson random measure on $[0,\infty)\times\N^*\times [0,1]$ of intensity measure:
\[
  \nd(\rmd s,\rmd i,\rmd \theta)
  \eqdef
  \lambdadm\,\rmd s\,\rmd i\,\rmd \theta
\]
(\/$\rmd i$ is the counting measure on $\N^*$, $\rmd s$ and $\rmd\theta$ are the Lebesgues measures on $[0,\infty)$ and $[0,1]$), and $\tNd=\Nd-\nd$ is the compensated measure.

\item
$\Nb$ is the Poisson random measure on $[0,\infty)\times\N^*\times\XX\times [0,1]$ of intensity measure:
\[
   \nb(\rmd s,\rmd i,\rmd z,\rmd \theta)
   \eqdef
   \lambdabm\,\kappa\,\bar D(z)\,\rmd s\,\rmd i\,\rmd z\, \rmd \theta\,,
\]
and  $\tNb=\Nb-\nb$ is the compensated measure.

\item
$\Nc$ is the Poisson random measure on $[0,\infty)\times\N^*\times\N^*\times [0,1]$ of intensity measure:
\[
   \nc(\rmd s,\rmd i,\rmd j,\rmd \theta)
   \eqdef
   \umax \,\rmd s\,\rmd i\,\rmd j\,\rmd \theta\,,
\]
and $\tNc=\Nc-\nc$ is the compensated measure.

\item
$(\mathbf{B}^i_t)_{t\geq 0}$, $i\geq 1$, is a family of mutually independent standard Brownian motions in $\R$.
\end{enumerate}
These four stochastic processes are mutually independent.

\end{enumerate}

\end{lemma}

\proof As the four basic mechanisms (natural death, birth, competition death, growth) are independent, we can write:
\begin{align}
\label{eq.lemma.Phi.nu.1}
  \Phi(\nu_t)
  =
  \Phi(\nu_0)
  +
  \Gd_t+\Gb_t+\Gc_t+\Gg_t
\end{align}
where $\Gd_t$, $\Gb_t$, $\Gc_t$, $\Gg_t$ are the terms associated with the natural death, the birth, the death due to competition, and the growth. We consider the four terms successively.

\subsubsection*{\it Death component $\Gd_t$}

The Monte Carlo procedure described in Section \ref{sec.monte.carlo} implies that:
\begin{align*}
  \Gd_t
  =
  \int_0^t\int_{\N^*}\int_0^1 
     \indic_{(i\leq N_{s^-})}\,
     \indic_{(\theta \leq \lambdad(x^i_{s^-})/\lambdadm)}\,
     [\Phi(\nu_{s^-}-\delta_{x^i_{s^-}}) - \Phi(\nu_{s^-})]\,
     \Nd(\rmd s,\rmd i,\rmd \theta)
\end{align*}
where $\Nd$ is a Poisson random measure described in \fenumii-\textit{a}. Here we use the labeling convention given at Remark \ref{rem.conv.numerotation}, $x^i_{s^-}$ is the $i^\textrm{th}$ particle of $\nu_{s^-}$. By introducing the compensated measure $\tNd=\Nd-\nd$ and the martingale $\Md_{\Phi,t}(\nu)$ defined by (\ref{eq.Md}), we get:
\begin{align*}
  \Gd_t
  &=
  \int_0^t \int_{\N^*} \int_0^1
     \indic_{(i\leq N_{s^-})}\,
     \indic_{(\theta \leq \lambdad(x^i_{s^-})/\lambdadm)}
     \,
     [\Phi(\nu_{s^-}-\delta_{x^i_{s^-}}) - \Phi(\nu_{s^-})]
     \;
     \lambdadm\,\rmd s\,\rmd i\,\rmd \theta
  +\Md_{\Phi,t}(\nu)
\\
  &=
  \int_0^t \int_{\N^*}
     \indic_{(i\leq N_{s^-})}\,
     \lambdad(x^i_{s^-})\,
     [\Phi(\nu_{s^-}-\delta_{x^i_{s^-}}) - \Phi(\nu_{s^-})]
     \;
     \rmd s\,\rmd i
  +\Md_{\Phi,t}(\nu)
\\
  &=
  \int_0^t \sum_{i=1}^{N_{s^-}} \lambdad(x^i_{s^-})\, 
  [\Phi(\nu_{s^-}-\delta_{x^i_{s^-}}) - \Phi(\nu_{s^-})]
  \;\rmd s
  +\Md_{\Phi,t}(\nu)
\\
  &=
  \int_0^t \int_{\XX} \lambdad(x)\,
     [\Phi(\nu_{s^-}-\delta_{x}) - \Phi(\nu_{s^-})]\,
    \nu_{s^-}(\rmd x)\;\rmd s
  +\Md_{\Phi,t}(\nu)
\end{align*}

\subsubsection*{\it Birth component $\Gb_t$}

The Monte Carlo procedure described in Section \ref{sec.monte.carlo} implies that:
\begin{align*}
  \Gb_t
  &=
  \int_0^t\int_{\N^*}\int_{\XX}\int_0^1 
     \indic_{(i\leq N_{s^-})}\,
     \indic_{\{\theta \leq (\lambdab(x^i_{s^-})\,D(x^i_{s^-},z))/(\lambdabm\,\kappa\,\bar D(z))\}}
     \,
\\
  &\qquad\qquad\qquad \qquad \qquad
    \times
     [\Phi(\nu_{s^-}+\delta_{(p_{x^i_{s^-}}+p_z,r_z)}) - \Phi(\nu_{s^-})]\,
     \Nb(\rmd s,\rmd i,\rmd z,\rmd \theta)
\end{align*}
where $\Nd$ is the Poisson random measure described in \fenumii-\textit{b}. 
By introducing the compensated measure $\tNb=\Nb-\nb$ and the martingale $\Mb_{\Phi,t}(\nu)$ defined by (\ref{eq.Mb}), we get:
\begin{align*}
  \Gb_t
  &=
  \int_0^t\int_{\N^*}\int_{\XX}\int_0^1 
     \indic_{(i\leq N_{s^-})}\,
     \indic_{\{\theta \leq (\lambdab(x^i_{s^-})\,D(x^i_{s^-},z))/(\lambdabm\,\kappa\,\bar D(z))\}}\,
\\
  &\qquad\qquad\qquad\qquad
    \times
     [\Phi(\nu_{s^-}+\delta_{(p_{x^i_{s^-}}+p_z,r_z)}) - \Phi(\nu_{s^-})]
     \;\lambdabm\,\kappa\,\bar D(z)\,\rmd s\,\rmd i\,\rmd z\, \rmd \theta
  +\Mb_{\Phi,t}(\nu)
\\
  &=
  \int_0^t
     \sum_{i=1}^{N_{s^-}}
       \lambdab(x^i_{s^-})\,
       \Big\{ \int_{\XX}  
                [\Phi(\nu_{s^-}+\delta_{(p_{x^i_{s^-}}+p_z,r_z)})-\Phi(\nu_{s^-})]
       \,D(x^i_{s^-},z)\,\rmd z\Big\}
     \;\rmd s
  +\Mb_{\Phi,t}(\nu)
\end{align*}

\subsubsection*{\it Competition component $\Gc_t$}

The Monte Carlo procedure described in Section \ref{sec.monte.carlo} implies that:
\begin{align*}
  \Gc_t
  &=
  \int_0^t\int_{\N^*}\int_{\N^*}\int_0^1 
     \indic_{(i\leq N_{s^-})}\,
     \indic_{(j\leq N_{s^-})}\,
     \indic_{(\theta \leq u(x^i_{s^-},x^j_{s^-})/\umax)}\,
\\
  &\qquad\qquad\qquad\qquad\qquad\qquad
    \times
     [\Phi(\nu_{s^-}-\delta_{x^i_{s^-}}) - \Phi(\nu_{s^-})]\,
     \Nc(\rmd s,\rmd i,\rmd j,\rmd \theta)
\end{align*}
where $\Nc$ is the Poisson random measure described in \fenumii-\textit{c}. 
By introducing the compensated measure $\tNc=\Nc-\nc$ and the martingale $\Mc_{\Phi,t}(\nu)$ defined by (\ref{eq.Mc}), we get:
\begin{align*}
  \Gc_t
  &=
  \int_0^t\int_{\N^*}\int_{\N^*}\int_0^1 
     \indic_{(i\leq N_{s^-})}\,
     \indic_{(j\leq N_{s^-})}\,
     \indic_{(\theta \leq u(x^i_{s^-},x^j_{s^-})/\umax)}\,
\\
  &\qquad\qquad\qquad\qquad
    \times
     [\Phi(\nu_{s^-}-\delta_{x^i_{s^-}}) - \Phi(\nu_{s^-})]\,
     \umax\,\rmd s\,\rmd i\,\rmd j\,\rmd \theta
  +\Mc_{\Phi,t}(\nu)
\\
  &=
  \int_0^t\int_{\N^*}\int_{\N^*} 
     \indic_{(i\leq N_{s^-})}\,
     \indic_{(j\leq N_{s^-})}\,
     u(x^i_{s^-},x^j_{s^-})\,
\\
  &\qquad\qquad\qquad\qquad
    \times
     [\Phi(\nu_{s^-}-\delta_{x^i_{s^-}}) - \Phi(\nu_{s^-})]\,
     \rmd s\,\rmd i\,\rmd j
  +\Mc_{\Phi,t}(\nu)
\\
  &=
  \int_0^t
     \sum_{i=1}^{N_{s^-}}
     \big[\sum_{j=1}^{N_{s^-}}  u(x^i_{s^-},x^j_{s^-})\big]\,
     [\Phi(\nu_{s^-}-\delta_{x^i_{s^-}}) - \Phi(\nu_{s^-})]\,
     \rmd s
  +\Mc_{\Phi,t}(\nu)\,.
\end{align*}

\subsubsection*{\it Growth component $\Gc_t$}

Consider now any instant between two instants of jumps (i.e. between two discrete events), for simplicity's sake we consider a time interval of $t\in[0,t_0]$ between 0 and the first instant of jump. For $t\in[0,t_0]$, the evolution of $\nu_t$ is modeled by the SDE (\ref{eq.growth}), that is:
\begin{align*}
  \rmd x^i_t = g(x^i_t,\nu_t)\,\rmd t
	               +\sigma(x^i_t,\nu_t)\,\rmd \mathbf{B}^i_t\,.
\end{align*}
For all $\Phi=(F,f)\in\DD$, from the Ito formula:
\begin{align*}
   f(x^i_t)
   &=
   f(x^i_0)
   + \int_0^t \GG f(x^i_s,\nu_s)\,\rmd s
   + \int_0^t \nabla f(x^i_s)\cdot \sigma(x^i_s,\nu_s)\,\rmd \mathbf{B}^i_s
\end{align*}
with the summation convention for repeated indices and where $a \eqdef \sigma\,\sigma^*$. Hence:
\begin{align*}
   \crochet{\nu_t,f}
   &=
   \crochet{\nu_0,f}
   + \int_0^t 
        \textstyle
        \crochet{
          \nu_s
          ,
          \GG f(\cdot,\nu_s)         }\,\rmd s
    +Y_t
\end{align*}
where
\begin{align*}
   Y_t
   &\eqdef
   \int_0^t \sum_{i=1}^{N_s}\nabla f(x^i_s)\cdot \sigma(x^i_s,\nu_s)\,\rmd \mathbf{B}^i_s
\end{align*}
This last expression is a martingale with quadratic variation:
\begin{align*}
   \crochet{Y}_t
   &\eqdef
   \int_0^t \sum_{i=1}^{N_s}\big|\nabla f(x^i_s)\cdot \sigma(x^i_s,\nu_s)\big|^2
       \,\rmd s
   =
   \int_0^t 
       \crochet{
           \nu_s
           ,
           |\nabla f\cdot \sigma(\cdot,\nu_s)|^2
       }
       \,\rmd s
\end{align*}
We get:
\begin{align*}
  F(\crochet{\nu_t,f})
  &=
  F(\crochet{\nu_0,f})
  +
  \int_0^t 
        \textstyle
        F'(\crochet{\nu_s,f})
        \times
        \crochet{
          \nu_s
          ,
          \GG f(\cdot,\nu_s)
          }\,\rmd s
\\
  &\qquad\qquad
  +
  \int_0^t 
        \textstyle
        F'(\crochet{\nu_s,f})\times
        \sum_{i=1}^{N_s}\nabla f(x^i_s)\cdot \sigma(x^i_s,\nu_s)\,\rmd \mathbf{B}^i_s
\\
  &\qquad\qquad
  +
  \demi
  \int_0^t 
        \textstyle
        F''(\crochet{\nu_s,f})
        \times
       \crochet{
           \nu_s
           ,
           |\nabla f\cdot \sigma(\cdot,\nu_s)|^2
       }
       \,\rmd s
\end{align*}
The infinitesimal generator associated with the growth phenomenon is obtained by taking the expectation of the previous expression. Finally
\begin{align*}
  \Phi(\nu_t)
  &=
  \Phi(\nu_0)
  +
  \int_0^t \LLg\Phi(\nu_{s^-})\,\rmd s
  +
  \Mg_{\Phi,t}(\nu)
\end{align*}
where $\Mg_{\Phi,t}(\nu)$ is the martingale defined in (\ref{eq.Mg}).

\medskip

Summing up these results in Equation (\ref{eq.lemma.Phi.nu.1}) ends the proof of the lemma.
\carre

\bigskip

As a corollary we now compute the predictable quadratic variation of the martingale processes $\Mtria_{f,t}(\nu)$:
\begin{corollary}
The predictable quadratic variation associated with martingale processes $\Md_.(\Phi)$, $\Mb_.(\Phi)$, $\Mc_.(\Phi)$, $\Mg_.(\Phi)$ defined in Equations (\ref{eq.M}) are:
\begin{subequations}
\label{eq.M.pqv}
\begin{align}
\label{eq.Md.pqv}
  \crochet{\Md_{\Phi,.}(\nu)}_t
  & =
  \int_0^t\int_{\XX}
     \lambdad(x)\,
       [\Phi(\nu_{s^-}-\delta_{x}) - \Phi(\nu_{s^-})]^2\,
  		\nu_{s^-}(\rmd x)\,\rmd s\,,
\\
\label{eq.Mb.pqv}
  \crochet{\Mb_{\Phi,.}(\nu)}_t
  & =
  \int_0^t\int_{\XX} 
     \lambdab(x)\,
     \int_{\XX} [\Phi(\nu_{s^-}+\delta_{(p_x+p_z,r_z)}) - \Phi(\nu_{s^-})]^2\,
     D(x,z)\,\rmd z\,
	  \,\nu_{s^-}(\rmd x)\, \rmd s\,,
\\
\label{eq.Mc.pqv}
  \crochet{\Mc_{\Phi,.}(\nu)}_t
  & =
  \int_0^t\int_{\XX}
    \Big[\int_{\XX}
     u(x,y)\,\nu_{s^-}(\rmd y)\Big]\,
     [\Phi(\nu_{s^-}-\delta_{x}) - \Phi(\nu_{s^-})]^2\,
   \nu_{s^-}(\rmd x)\,\rmd s\,,
\\
\label{eq.Mg.pqv}
  \crochet{\Mg_{\Phi,.}(\nu)}_t
  &=
  \int_0^t 
  \int_\XX
      \big|F'(\crochet{\nu_{s^-},f})\big|^2
      \,
      \big| \nabla f(x)\cdot \sigma(x,\nu_{s^-})\big|^2
      \,
      \nu_{s^-}(\rmd x)\,\rmd s
        \,.
\end{align}
\end{subequations}
and by independence of the processes $\tNd$, $\tNb$, $\tNc$ and $\mathbf{B}^i$, $i\geq 1$, we get:
\[
\crochet{\mathbf{M}_{\Phi,.}(\nu)}_t
=
\crochet{\Md_{\Phi,.}(\nu)}_t
+\crochet{\Mb_{\Phi,.}(\nu)}_t
+\crochet{\Mc_{\Phi,.}(\nu)}_t
+\crochet{\Mg_{\Phi,.}(\nu)}_t
\,.
\]
\end{corollary}

\proof
Consider a martingale process:
\[
  \mathbf{M}_t \eqdef
  \int_0^t\int_E \kappa(s,x)\,\mathbf{\tilde N}(\rmd s,\rmd x)
\]
where $\mathbf{\tilde N}$ is the compensated measure associated with a Poisson random measure $\mathbf{N}$ on $[0,\infty)\times E$ of intensity measure $\mathbf{n}$ (i.e. $\mathbf{\tilde N}= \mathbf{N}-\mathbf{n}$), and $\kappa(s,x)$ is a predictable process such that $\E\int_0^T\int_E \kappa^2(s,x)\,\mathbf{n}(\rmd s,\rmd x)<\infty$. Then:
\[
  \crochet{\mathbf{M}}_t
  =
  \int_0^t\int_E \kappa^2(s,x)\,\mathbf{n}(\rmd s,\rmd x)\,.
\]
(see \cite[Section II-3]{ikeda1981a}). We apply this result for example for the first term, from (\ref{eq.Md}):
\begin{align*}
  \crochet{\Md_{\Phi,.}(\nu)}_t
  & =
  \int_0^t\int_{\N^*}\int_0^1 
     \indic_{(i\leq N_{s^-})}\,
     \indic_{(\theta \leq \lambdad(x^i_{s^-})/\lambdadm)}\,
\\
  &\qquad\qquad \qquad\qquad
    \times
       [\Phi(\nu_{s^-}-\delta_{x^i_{s^-}}) - \Phi(\nu_{s^-})]^2\;
     \nd(\rmd s,\rmd i,\rmd \theta)
\\ 
  & =
  \int_0^t\int_{\N^*}\int_0^1 
     \indic_{(i\leq N_{s^-})}\,
     \indic_{(\theta \leq \lambdad(x^i_{s^-})/\lambdadm)}\,
       [\Phi(\nu_{s^-}-\delta_{x^i_{s^-}}) - \Phi(\nu_{s^-})]^2\,
     \lambdadm\,\rmd \theta\,\rmd i\,\rmd s\,,
\end{align*}
integrating in $\theta$ leads to
\begin{align*}
  \crochet{\Md_{\Phi,.}(\nu)}_t
  & =
  \int_0^t\int_{\N^*}
     \textstyle
     \indic_{(i\leq N_{s^-})}\,
     \frac{\lambdad(x^i_{s^-})}{\lambdadm}\,
       [\Phi(\nu_{s^-}-\delta_{x^i_{s^-}}) - \Phi(\nu_{s^-})]^2\,
     \lambdadm\,\rmd i\,\rmd s
\\
  & =
  \int_0^t\int_{\N^*}
     \textstyle
     \indic_{(i\leq N_{s^-})}\,
     \lambdad(x^i_{s^-})\,
       [\Phi(\nu_{s^-}-\delta_{x^i_{s^-}}) - \Phi(\nu_{s^-})]^2\,\rmd i\,
     \rmd s
\\
  & =
  \int_0^t\int_{\XX}
     \lambdad(x)\,
       [\Phi(\nu_{s^-}-\delta_{x}) - \Phi(\nu_{s^-})]^2\,
     \nu_{s^-}(\rmd x)\,\rmd s
\end{align*}
The same approach could be applied to the next two expressions, the last assertion (\ref{eq.Mg.pqv}) is due to the properties of the Brownian motion.
\carre

\bigskip

We now consider a particular case of test functions $\Phi$:
\[
   \Phi(\nu) = \crochet{\nu,f}
\]
(i.e. with $F=$id.) for any function $f:\XX\mapsto\R$ twice continuously differentiable, bounded with bounded derivatives.

We define:
\begin{align*}
   \lltria f(\nu)
   &\eqdef
   \LLtria \Phi(\nu)
\\
  \mtria_{f,t}(\nu)
  &\eqdef
  \Mtria_{\Phi,t}(\nu)
\end{align*}
for $\vartriangle=\textrm{\rm d,b,c,g}$. Hence:
\begin{align}
\nonumber
  \crochet{\nu_t,f}
  &=
  \crochet{\nu_0,f}
  +
  \int_0^t \ell f(\nu_{s})\;\rmd s+ \mathbf{m}_{f,t}(\nu)
\\
\label{eq.nu.f.total}
  &=
  \crochet{\nu_0,f}
  +
  \sum_{\vartriangle=\textrm{\rm d,b,c,g}}
    \Big\{
       \int_0^t \lltria f(\nu_{s})\;\rmd s
       +
       \mtria_{f,t}(\nu)
    \Big\}
\end{align}
where
\begin{enumerate}

\item
$\ell=\lld+\llb+\llc+\llg$ are the infinitesimal generators defined by:
\begin{align*}
  \lld f(\nu)
  & =
  -\int_{\XX} \lambdad(x)\,f(x)\,\nu(\rmd x)\,,
\\
  \llb f(\nu)
  &=
  \int_\XX
       \textstyle
       \lambdab(x)\,
       \big[ \int_{\XX}  f((p_x+p_z,r_z))\, D(x,z)\,\rmd z\big]
   \,\nu(\rmd x)\,,
\\
  \llc f(\nu)
  &=
  -\int_\XX 
       \textstyle
       \big[\int_\XX u(x,y)\,\,\nu(\rmd y)\big] \,
       f(x)\,\,\nu(\rmd x)\,,
\\
  \llg f(\nu)
  &=
  \int_\XX \GG f(x,\nu)    \,\nu(\rmd x)
\end{align*}
corresponding respectively to natural death, birth, death by competition and growth respectively. 

\item
the martingale terms $\mathbf{m}_{f,t}(\nu)=\md_{f,t}(\nu)+\mb_{f,t}(\nu)+\mc_{f,t}(\nu)+\mg_{f,t}(\nu)$ are defined by:
\begin{align*}
  \md_{f,t}(\nu)
  & =
  -\int_0^t\int_{\N^*}\int_0^1 
     \indic_{(i\leq N_{s^-})}\,
     \indic_{(\theta \leq \lambdad(x^i_{s^-})/\lambdadm)}\, f(x^i_{s^-})\,
     \tNd(\rmd s,\rmd i,\rmd \theta)\,,
\\
  \mb_{f,t}(\nu)
  & =
  \int_0^t\int_{\N^*}\int_{\XX}\int_0^1 
     \indic_{(i\leq N_{s^-})}\,
     \indic_{(\theta \leq \lambdab(x^i_{s^-})/\lambdabm)}\,
     f(x^i_{s^-}+z)\,
     \tNb(\rmd s,\rmd i,\rmd z,\rmd \theta)\,,
\\
  \mc_{f,t}(\nu)
  & =
  -\int_0^t\int_{\N^*}\int_{\N^*}\int_0^1 
     \indic_{(i\leq N_{s^-})}\,
     \indic_{(j\leq N_{s^-})}\,
     \indic_{(\theta \leq u(x^i_{s^-},x^j_{s^-})/\umax)}\,
\\
  &\qquad\qquad\qquad\qquad\qquad\qquad
    \times
     f(x^i_{s^-})\,
     \tNc(\rmd s,\rmd i,\rmd j,\rmd \theta)
     \,,
\\
  \mg_{f,t}(\nu)
  &=
  \int_0^t 
        \textstyle
        \sum_{i=1}^{N_s}\nabla f(x^i_s)\cdot \sigma(x^i_s,\nu_s)\,\rmd \mathbf{B}^i_s
        \,,
\end{align*}
where the compensated Poisson random measures $\tNd$, $\tNb$, $\tNc$ and the Brownian motions $(\mathbf{B}^i_t)_{t\geq 0}$, $i\geq 1$, are introduced in Lemma \ref{lemma.Phi.nu}.
The corresponding predictable quadratic variation terms are:
\begin{subequations}
\label{eq.m.pqv}
\begin{align}
\label{eq.md.pqv}
  \crochet{\md_{f,.}(\nu)}_t
  & =
  \int_0^t\int_{\XX}
     \lambdad(x)\,
       f^2(x)\,
  		\nu_{s}(\rmd x)\,\rmd s\,,
\\
\label{eq.mb.pqv}
  \crochet{\mb_{f,.}(\nu)}_t
  & =
  \int_0^t\int_{\XX} 
     \textstyle
     \lambdab(x)\,
       \big[\int_{\XX} f^2((p_x+p_z,r_z)) \, D(x,z)\,\rmd z\big]
     \,\,\nu_{s}(\rmd x)\, \rmd s\,,
\\
\label{eq.mc.pqv}
  \crochet{\mc_{f,.}(\nu)}_t
  & =
  \int_0^t\int_{\XX}
     \textstyle
    \big[\int_{\XX} u(x,y)\,\nu_{s}(\rmd y)\big]\, f^2(x)\,
   \nu_{s}(\rmd x)\,\rmd s\,,
\\
\label{eq.mg.pqv}
  \crochet{\mg_{f,.}(\nu)}_t
  &=
  \int_0^t 
        \int_\XX
          \big|\nabla f(x)\cdot \sigma(x,\nu_s)\big|^2\,\nu_s(\rmd x)
        \,\rmd s
        \,.
\end{align}
\end{subequations}
and by independence of the processes $\tNd$, $\tNb$, $\tNc$ and $\mathbf{B}^i$, $i=\geq 1$, we get:
\[
\crochet{\mathbf{m}_{f,.}(\nu)}_t
  =
  \crochet{\md_{f,.}(\nu)}_t
  +\crochet{\mb_{f,.}(\nu)}_t
  +\crochet{\mc_{f,.}(\nu)}_t
  +\crochet{\mg_{f,.}(\nu)}_t\,.
\]
\end{enumerate}

\subsection{Control of the size of the population}

\begin{lemma}
\label{lemma.pop.finie}
If for some $p\geq 1$, $\E (\crochet{\nu_0,1}^p)<\infty$, then for any $T>0$,
\[
  \E \sup_{0\leq t\leq T} \crochet{\nu_t,1}^p<\infty\,.
\]
\end{lemma}

\proof
We apply Lemma \ref{lemma.Phi.nu} with $F(\xi)=\xi^p$ and $f(x)=1$.
As the terms corresponding to death (natural death, competition) are less than zero, and as $\LLg\Phi=0$, we get:
\begin{align*}
  N_t^p
  &\leq
  N_0^p
  +
  \int_0^t\int_{\N^*}\int_{\XX}\int_0^1 
     \indic_{(i\leq N_{s^-})}\,
     \indic_{(\theta \leq \lambdab(x^i_{s^-})/\lambdabm)}\,
     [(N_{s^-}+1)^p-N_{s^-}^p]\,
     \Nb(\rmd s,\rmd i,\rmd z,\rmd \theta)\,.
\end{align*}
We introduce the stopping time:
\[
  \tau_n \eqdef \inf\{t\geq 0\,;\, N_t\geq n\}
\]
then
\begin{align*}
  \sup_{0\leq t\leq T\wedge \tau_n}N_t^p
  &\leq
  N_0^p
  +
  \int_0^{T\wedge \tau_n}\int_{\N^*}\int_{\XX}\int_0^1 
     \indic_{(i\leq N_{s^-})}\,
     \indic_{(\theta \leq \lambdab(x^i_{s^-})/\lambdabm)}\,
     [(N_{s^-}+1)^p-N_{s^-}^p]\,
\\
  &\qquad\qquad\qquad\qquad\qquad\qquad\qquad\qquad\qquad
     \times  \Nb(\rmd s,\rmd i,\rmd z,\rmd \theta)
\end{align*}
as $(k+1)^p-k^p\leq C_p\,(1+k^{p-1})$ for some constant $C_p$, we get:
\begin{align*}
  \sup_{0\leq t\leq T\wedge \tau_n}N_t^p
  &\leq
  N_0^p
  +
  C_p \int_0^{T\wedge \tau_n}\int_{\N^*}\int_{\XX}\int_0^1 
     \indic_{(i\leq N_{s^-})}\,
     [1+N_{s^-}^{p-1}]\,
  \Nb(\rmd s,\rmd i,\rmd z,\rmd \theta)\,.
\end{align*}
Taking expectation leads to:
\begin{align}
\nonumber
  \E\sup_{0\leq t\leq T\wedge \tau_n}N_t^p
  &\leq
  \E (N_0^p)
  +
  C_p \;\E \int_0^{T\wedge \tau_n}\int_{\N^*}\int_{\XX}\int_0^1 
     \indic_{(i\leq N_{s^-})}\,
     [1+N_{s^-}^{p-1}]\,
     \nb(\rmd s,\rmd i,\rmd z,\rmd \theta)
\\
\nonumber
  &\leq
  \E (N_0^p)
  +
  C_p \int_0^T 
  \E\big(
       \indic_{(s\leq T\wedge \tau_n)}\,N_s\,[1+N_s^{p-1}] 
  \big)\, \rmd s
\\
\nonumber
  &\leq
  \E (N_0^p)
  +
  C_p \int_0^T 
  \E\big(
       \indic_{(s\leq T\wedge \tau_n)}\, [1+N_s^p] \big)\, \rmd s
\\
\nonumber
  &\leq
  \E (N_0^p)
  +
  2\,C_p \int_0^T 
     \E(
       N_{s\wedge \tau_n}^p)\, \rmd s
\\
\label{eq.lemma.pop.finie.1}
  &\leq
  \E (N_0^p)
  +
  2\,C_p \int_0^T 
  \E\big(
       \sup_{0\leq s\leq t\wedge \tau_n}N_{s\wedge \tau_n}^p
  \big)\, \rmd t
\end{align}
and by Gronwall's lemma
\begin{align}
\label{eq.lemma.pop.finie}
  \E\sup_{0\leq t\leq T\wedge \tau_n}N_t^p
  &\leq C_{p,T}
\end{align}
for all $n$. We new want to check that $\tau_n\to\infty$ a.s. Suppose that $\tau_n\not\to \infty$, then there exists $T_0<\infty$ such that $\epsilon_0=\P(\sup_n \tau_n \leq T_0)>0$, then:
\[
  \E\big(\sup_{0\leq t\leq T_0\wedge \tau_n}N_t^p\big)
  \geq 
  n^p\,\epsilon_0
\]
which is a contradiction. Letting $n\to\infty$ in (\ref{eq.lemma.pop.finie}), by Fatou's lemma we prove the result.

Note that in Equation (\ref{eq.lemma.pop.finie.1}) the different constants $C_p$ depend only on $\lambdabm$ and on universal constants, but do not depend on the functions $u$, $g$, $\sigma$.
\carre

\section{Large population limit}
\label{sec.asymptotic}

Let $k$ be the initial population size, i.e. $k=\crochet{\nu_0,1}$,
and replace $u$ by $u^k$, $g$ by $g^k$ and $\sigma$ by
$\sigma^k$. Let $(\nu^k_t)_{0\leq t\leq T}$ be the Markov process defined in the
previous section with initial population size $k$. We define: 
\[
   \mu^k_t \eqdef \frac{1}{k}\,\nu^k_t\,.
\]
In this section we study the asymptotic property of the law of the process $(\mu^k_t)_{0\leq t\leq T}$ on the space $\D([0,T],\MM_F(\XX))$ of \emph{càdlàg} functions from $[0,T]$ with values in $\MM_F(\XX)$ (\footnote{If not mentioned $\D([0,T],\MM_F(\XX))$ is equipped with the Skorohod topology associated with the weak topology on $\MM_F(\XX)$ see appendix.}). According to the general approach depicted in the appendix, a classical method to do this consists of deducing both convergence and characterization of the limiting process from the convergence and properties of $ \crochet{\mu^k_t,f}$, for all $f$ in a suitable class. 

\medskip

From (\ref{eq.nu.f.total}):
\begin{align}
\label{eq.nu.k.f.total.1}
  \crochet{\nu^k_t,f}
  &=
  \crochet{\nu^k_0,f}
  +
  \int_0^t \ell^k f(\nu^k_{s})\;\rmd s
  +
  \mathbf{m}_{f,t}^k(\nu^k)
\end{align}
where:
\begin{align*}
  \ell^k f(\nu)
  & =
  \int_{\XX} 
  \textstyle
  \Big\{ - \lambdad(x)\,f(x)
       + \lambdab(x)\, \big[ \int_{\XX}  f((p_x+p_z,r_z))\,D(x,z)\,\rmd z\big]
\\[-0.2em]
  &\qquad\qquad\qquad
    \textstyle
       - \big[\int_\XX u^k(x,y)\,\,\nu(\rmd y)\big] \,f(x)
    +
         \GG^k f(x,\nu)
    \Big\}
    \;\nu(\rmd x) 
    \vphantom{\displaystyle\int_{\XX}}
\end{align*}
with
\begin{align*}
  \GG^k f(x,\nu)
  \eqdef\textstyle
  \nabla f(x)\cdot g^k(x,\nu)
         + \demi\,\frac{\partial^2 f}{\partial x_\ell\,\partial x_{\ell'}}
             \,a^k_{\ell\ell'}(x,\nu)
\end{align*}
and
\begin{align*}
  &
  \mathbf{m}^k_{f,t}(\nu)
  =
  -\int_0^t\int_{\N^*}\int_0^1 
     \indic_{(i\leq N_{s^-})}\,
     \indic_{(\theta \leq \lambdad(x^i_{s^-})/\lambdadm)}\, f(x^i_{s^-})\,
     \tNd(\rmd s,\rmd i,\rmd \theta)
\\
  & \qquad +
  \int_0^t\int_{\N^*}\int_{\XX}\int_0^1 
     \indic_{(i\leq N_{s^-})}\,
     \indic_{(\theta \leq \lambdab(x^i_{s^-})/\lambdabm)}\,
     f((p_{x^i_{s^-}}+p_z,r_z))\,
     \tNb(\rmd s,\rmd i,\rmd z,\rmd \theta)
\\
  & \qquad  
  -\int_0^t\int_{\N^*}\int_{\N^*}\int_0^1 
     \indic_{(i\leq N_{s^-})}\,
     \indic_{(j\leq N_{s^-})}\,
     \indic_{(\theta \leq u^k(x^i_{s^-},x^j_{s^-})/\umax)}\,
     f(x^i_{s^-})\,
     \tNc(\rmd s,\rmd i,\rmd j,\rmd \theta)
\\
  & \qquad +
  \int_0^t 
        \textstyle
        \sum_{i=1}^{N_s}\nabla f(x^i_s)\cdot \sigma^k(x^i_s,\nu_s)\,\rmd \mathbf{B}^i_s
\end{align*}
Dividing (\ref{eq.nu.k.f.total.1}) by $k$ leads to:
\begin{align}
\nonumber
  \crochet{\mu^k_t,f}
  &=
  \crochet{\mu^k_0,f}
  +
  \int_0^t \frac{1}{k}\,\ell^k f(\nu^k_{s})\;\rmd s
  +
  \frac{1}{k}\,\mathbf{m}_{f,t}^k(\nu^k)
\\
\label{eq.crochet.mu.k}
  &=
  \crochet{\mu^k_0,f}
  +
  \int_0^t \frac{1}{k}\,\ell^k f(k\,\mu^k_{s})\;\rmd s
  +
  \frac{1}{k}\,\mathbf{m}_{f,t}^k(k\,\mu^k)
\end{align}
Define:
\begin{align}
\label{eq.Z.k}
   Z^k_t \eqdef \frac{1}{k}\,\mathbf{m}_{f,t}^k(k\,\mu^k)
\end{align}
and 
\begin{subequations}
\label{eq.bll-total} 
\begin{align}
\label{eq.bll} 
 \bll &= \blld+ \bllb+ \bllc+ \bllg
\end{align}
with
\begin{align}
\label{eq.blld} 
  \blld f(\mu) 
  &= 
  \frac{1}{k} \lld f(k\, \mu) = \lld f(\mu) 
  =
  - \int_{\XX} \lambdad(x)\,f(x) \,\mu(\rmd x)\,,
\\
\label{eq.bllb} 
  \bllb f(\mu) 
  &= 
  \frac{1}{k} \llb f(k\, \mu) = \llb f(\mu) 
  =
  \int_{\XX} 
       \lambdab(x)\,\textstyle
       \big[ \int_{\XX}  f((p_x+p_z,r_z))\,D(x,z)\,\rmd z\big]
    \;\mu(\rmd x)\,,
\\
\label{eq.bllc} 
  \bllc f(\mu) 
  &= 
  -\int_\XX \textstyle
      \big[\int_\XX \bar u(x,y)\,\mu(\rmd y)\big] \, f(x)\,\mu(\rmd x) 
\\
\label{eq.bllg} 
  \bllg f(\mu)  
  &=
  \int_\XX 
	   \bar\GG f(x,\mu) \,\mu(\rmd x)
\end{align} 
\end{subequations}
and
\begin{align*}
  \bar\GG f(x,\mu) 
  &\eqdef
	   \nabla f(x)\cdot \bar g(x,\mu) 
	   +
	   \textstyle\demi\,\frac{\partial^2 f(x)}{\partial x_\ell\,\partial x_{\ell'}}
	             \,\bar a_{\ell\ell'}(x,\mu)\,.
\end{align*} 
Hence (\ref{eq.crochet.mu.k}) reads:
\begin{align}
\label{eq.mu.k}
  \crochet{\mu^k_t,f}
  &=
  \crochet{\mu^k_0,f}
  +
  \int_0^t \bll f(\mu^k_{s})\;\rmd s
  +
  Z^k_t
  +
  R^k_t
\end{align}
where
\begin{align}
\label{eq.R.k}
  R^k_t
  \eqdef
  \int_0^t 
    \Big[
      \frac{1}{k}\,\ell^k f(k\,\mu^k_{s})
      -
      \bll f(\mu^k_{s})
     \Big]\,\rmd s\,.
\end{align}
From (\ref{eq.m.pqv}) we get:
\begin{align}
\nonumber
  \crochet{Z^k}_t
  &=
  \frac{1}{k^2}\,\Crochet{\mathbf{m}_{f,.}^k(k\,\mu^k)}_t
\\
\nonumber
  &=
  \frac{1}{k}  
  \int_0^t \int_\XX \Big\{ 
    \textstyle
    \lambdad(x)\,f^2(x)\,
	+
	\lambdab(x) \, \big[\int_\XX f^2((p_x+p_z,r_z))\, D(x,z) \, \rmd z\big] 
\\
\nonumber
  & \qquad\qquad\qquad\qquad
    \textstyle
	+ \big[\int_\XX k\, u^k(x,y) \,  \mu^k_{s}(\rmd y)\big]\,f^2(x)
\\
\label{eq.crochet.Z.k}
  & \qquad\qquad\qquad\qquad
    \textstyle
	+ |\nabla f(x)\cdot \sigma ^k(x,k\, \mu_s^k)|^2
	\Big\}  
	\, \mu^k_{s}(\rmd x)\, \rmd s \,.
\end{align}

We are now ready to state the main result of this section.
\begin{theorem}
\label{th.cv}
Suppose that:
\begin{subequations}
\label{hyp.cv}
\begin{align}
\label{hyp.cv.u}
   k\,u^k(x,y)        
   &
   \xrightarrow[k\to \infty]{L^\infty(\XX\times\XX)} 
   \bar u(x,y)\,,
   \\
\label{hyp.cv.g}
   g^k(x,k\,\mu)      
   &
   \xrightarrow[k\to \infty]{L^\infty(\XX\times\MM_F(\XX))} 
   \bar g(x,\mu)\,,
   \\
\label{hyp.cv.sigma}
   \sigma^k(x,k\,\mu) 
   &\xrightarrow[k\to \infty]{L^\infty(\XX\times\MM_F(\XX))} 
   \bar\sigma(x,\mu)
\end{align}
\end{subequations}
and that $\bar g(x,\xi)$ and $\bar a(x,\xi)$ are bounded and Lipschitz in $x$ uniformly in $\xi$, i.e.
\begin{align}
\label{eq.Lipschitz}
  |\bar g(x,\xi)-\bar g(x,\xi')|
  +
  |\bar a(x,\xi)-\bar a(x,\xi')|
  \leq 
  C\,\rho(\mu,\mu')
\end{align}
where $\rho$ is the Prohorov metric on $\MM_F(\XX)$ ($\rho$ generates the topology of weak convergence on  $\MM_F(\XX)$). Suppose also that:
\begin{align}
\label{eq.condi.ini}
  \mu^k_0
  \xrightarrow[k\to\infty]{\textrm{\rm law}}
  \xi_0\,.
\end{align}
where $\xi_0\in\MM_F(\XX)$ (deterministic). 

Then $(\mu ^k)_{k \in \bbN}$ converges in law to
a deterministic process $\xi \in  \CC([0,T],\MM_F(\XX))$, characterized
by
\begin{align}
\label{eq.xi.caract}
  \crochet{\xi_t,f}
  &=
  \crochet{\xi_0,f} + \int_0^t \bll f(\xi_s)\,\rmd s
\end{align}
where $\bll f$ is defined in Equations (\ref{eq.bll-total}). 
\end{theorem} 

\begin{remark}
Hypothesis (\ref{hyp.cv.u}) can be understood as a ``small competition/large population'' asymptotic. In the example of Section \ref{sec.foret},
the term $g^k$ is given by (\ref{modele.richards}) and the convergence (\ref{hyp.cv.g}) can be deduced from (\ref{hyp.cv.u}) as:
\[
  g^k(x,k\,\mu)
  =
  \Psi\big({\textstyle\int k\,u^k(x,y)\,\mu(\rmd y)}\big)\,  R(r)
  \xrightarrow[k\to\infty]{}
  \Psi\big({\textstyle\int \bar u^k(x,y)\,\mu(\rmd y)}\big)\,  R(r)
  \,.
\] 
Note also that in this example $\sigma^k(x,k\,\mu)\equiv 0$.
\end{remark}

The end of the section is devoted to the proof of Theorem \ref{th.cv}.

\subsubsection*{Uniqueness}

\begin{lemma}
Equation (\ref{eq.xi.caract}) admits a unique solution in $\CC([0,T],\MM_F(\XX))$.
\end{lemma}

\proof
First note that, as for the proof of Lemma \ref{lemma.pop.finie}:
\begin{align*}
  \crochet{\xi_t,1}
  &\leq
  \crochet{\xi_0,1}
  +\int_0^t \int_{\XX} 
       \lambdab(x)\,\xi_s(\rmd x)\,\rmd s
  \leq
  \crochet{\xi_0,1}
  +
  \lambdabm\, \int_0^t \crochet{\xi_s,1}\,\rmd s
\end{align*} 
so that, from Gronwall's lemma, $\sup_{t\in[0,T]}\crochet{\xi_t,1} \leq \crochet{\xi_0,1}\,\exp(\lambdabm\,T)=:\kappa_T$.

\medskip



Consider two solutions $\xi_t$ and $\tilde\xi_t$ of (\ref{eq.xi.caract}) with the same initial condition $\xi_0$, we prove that $\Normtv{\xi_t-\tilde\xi_t}=0$ for all $t\in[0,T]$ where:
\[
   \Normtv{\mu}
   \eqdef
   \sup_{B\in \BB(\XX)} |\mu(B)|
   =
   \sup_{\substack{f\textrm{ Borel function}\\ \Norm{f}_\infty\leq 1}} |\crochet{\mu,f}|
   =
   \sup_{\substack{f\in\CC^\infty_0(\XX) \\ \Norm{f}_\infty\leq 1}} |\crochet{\mu,f}|
\]
where $\CC^\infty_0(\XX)$ is the set of smooth functions with compact support. For any $f\in\CC^\infty_0(\XX)$ with $\Norm{f}_\infty\leq 1$, we have:
\begin{align}
\nonumber
  &|\crochet{\xi_t-\tilde\xi_t,f}|
  \leq
  \int_0^t\Big\{\Big|\int_{\XX} 
     \textstyle\lambdad(x)\,f(x) \, 
     [\xi_s(\rmd x)-\tilde\xi_s(\rmd x)]
  \Big|
\\
\nonumber
  &\qquad
  +
  \Big|\int_{\XX} 
       \textstyle\lambdab(x)\,\textstyle
       \big( \int_{\XX}  f((p_x+p_z,r_z))\,D(x,z)\,\rmd z\big)\, 
    \; [\xi_s(\rmd x)-\tilde\xi_s(\rmd x)]\Big|
\\
\nonumber
  &\qquad
  + \Big|\int_\XX \textstyle
      \big(\int_\XX \bar u(x,y)\, \xi_s(\rmd y)\big) \, f(x)\, 
      [\xi_s(\rmd x)-\tilde\xi_s(\rmd x)] \Big|
\\
\nonumber
  &\qquad
  + \Big|\int_\XX \textstyle
     \big(
         \int_\XX \bar u(x,y)\, [\xi_s(\rmd y)-\tilde\xi_s(\rmd y)]
      \big) \, f(x)
      \, \tilde\xi_s(\rmd x) \Big|
\\
\nonumber
  &\qquad
  +
  \Big|\int_\XX 
     \big(
	   \nabla f(x)\cdot \bar g(x, \xi_s) 
	   +
	   \textstyle\demi\,\frac{\partial^2 f(x)}{\partial x_\ell\,\partial x_{\ell'}}
	             \,\bar a_{\ell\ell'}(x, \xi_s)
	 \big) \, [\xi_s(\rmd x)-\tilde\xi_s(\rmd x)]\Big|
\\
\label{eq.xi-xitilde}
  &\qquad
  +
  \Big|\int_\XX 
     \big(
	   \nabla f(x)\cdot [\bar g(x, \xi_s)-\bar g(x, \tilde\xi_s)] 
	   +
	   \textstyle\demi\,\frac{\partial^2 f(x)}{\partial x_\ell\,\partial x_{\ell'}}
	             \,[\bar a_{\ell\ell'}(x, \xi_s)-\bar a_{\ell\ell'}(x, \tilde\xi_s)]
	 \big) \, \tilde\xi_s(\rmd x)\Big|\Big\}\,\rmd s
\end{align} 
The first five terms of the {\sc rhs} of (\ref{eq.xi-xitilde}) are controlled the same way, for example the third term is:
\begin{align*}
  &\umaxb\,\kappa_T\,
     \Big|\int_\XX \textstyle
       \underbrace{\textstyle\frac{1}{\umaxb\,\kappa_T}
      \,\big(\int_\XX \bar u(x,y)\, \xi_s(\rmd y)\big) \, f(x)}_{\norm{\cdot}_\infty\leq 1}\, 
      [\xi_s(\rmd x)-\tilde\xi_s(\rmd x)] \Big|
  \leq \umaxb\,\kappa_T\,\Normtv{\xi_s-\tilde\xi_s}
\end{align*} 
where $\bar u(x,y)\leq \umaxb$.
For the last term of the {\sc rhs} of (\ref{eq.xi-xitilde}) we use the fact that the coefficients $\bar g(x,\xi)$ and $\bar a(x,\xi)$ are bounded Lipschitz in $x$ uniformly in $\xi$. We get:
\begin{align*}
  &|\crochet{\xi_t-\tilde\xi_t,f}|
  \leq
  C_T\, \int_0^t \Normtv{\xi_s-\tilde\xi_s}\,\rmd s
\end{align*} 
so, thanks to  Grownwall's lemma, $\Normtv{\xi_t-\tilde\xi_t}=0$ for all $t\in[0,T]$.
\carre

\subsubsection*{Moment estimate}

There exist a constant $C_T$ depending only on $T$ such that:
\begin{align}
\label{eq.lemma.pop.finie.unif}
  \sup_{k\in\N} \E \Big(\sup_{0\leq t\leq T}\crochet{\mu^k_t, 1} ^3\Big) 
  &\leq C_T
\end{align}
This result is similar to Lemma \ref{lemma.pop.finie}. Indeed, by dividing each term of Equation (\ref{eq.lemma.pop.finie.1}) by $k$, we can apply the same reasoning thanks to the remark at the end of the proof of the lemma.

\subsubsection*{Tightness for fixed $f$} 

We first show that the family of laws of $\crochet{\mu^k_\cdot, f}$ is
tight, for each $f\in \CC^2_b(\XX)$.  For fixed $t<T$, since $f$ is bounded, we
have 
\[
  \P\big(|\crochet{\mu^k_t, f}| > K\big) 
  \leq 
  \frac{1}{K}\,C_f\, 
  \sup_{k} \E \Big(\sup_{0\leq t\leq T}  \crochet{\mu^k_t, 1} \Big)
\]
Using~(\ref{eq.lemma.pop.finie.unif}), we conclude that the sequence
of the laws of $(\crochet{\mu^k_t, f} )_{k\in\N}$ is tight. Then
denote by $A^k_t$ the finite variation part of $\crochet{\mu^k_t, f}$. We have
\[
  A^k_t
  =
  \sum_{\vartriangle=\textrm{\rm d,b,c,g}} 
  \int_0 ^t \frac{1}{k} \, \lltria f(k\, \mu^k_{s})\;\rmd s 
\]
with
\begin{align*}
  \Big| \frac{1}{k} \lld f(k\, \mu) \Big|
  &= 
  \Big|\int _\XX \lambdad (x)\, f(x) \, \mu(\rmd x) \Big|
  \leq C_f \, \crochet{\mu,1}
\\
  \Big| \frac{1}{k} \llb f(k\, \mu) \Big|
  &= 
  \Big|\int_\XX 
    \lambdab (x)\, 
    \textstyle
    \big[\int_\XX f((p_x+p_z,r_z)) \, D(x,z)\, \rmd z \big]\, 
	\mu(\rmd x) \Big|
  \leq C_f \, \crochet{\mu,1}
\\
  \Big|\frac{1}{k} \llc f(k\, \mu) \Big|
  &= 
  \Big|\int _\XX	
       \textstyle
		\big[\int_\XX k\, u^k(x,y)\, \mu(\rmd y)\big]\, 
		\mu(\rmd y) \Big| 
  \leq C_f \, \crochet{\mu,1}^2
\\
  \Big|\frac{1}{k} \llg f(k\, \mu) \Big|
  &=
  \Big|\int_\XX
     \big[
       \nabla f(x)\cdot g^k(x,k\, \mu) 
       +
       \textstyle\demi\,\frac{\partial^2 f(x)}{\partial x_\ell\,\partial x_{\ell'}}
             \,a_{\ell\ell'}^k (x,k\, \mu)
     \big]
  \,\mu(\rmd x)
  \Big|
  \leq  C_f \, \crochet{\mu,1} \ .
\end{align*} 
Hence for any sequence $(\tau_n)_{n\in\N}$ of stopping times bounded by $T$,
\[
  \E |A_{\tau_k + \theta}^k -  A_{\tau_k}^k | 
  \leq 
  C \, \theta 
\]
using~(\ref{eq.lemma.pop.finie.unif}) again. Similarly, we can prove
\[
 \E |\crochet{Z^k}_{\tau_k + \theta} - \crochet{Z^k}_{\tau_k} |
  \leq  
  \frac{C \, \theta}{k} 
\]
thanks to (\ref{eq.crochet.Z.k}) and (\ref{eq.lemma.pop.finie.unif}). According to the Aldous--Rebolledo criteria, this ensures the tightness of the laws of 
$\crochet{\mu_.^k,f}$.

\subsubsection*{Tightness}

The last result allows us to apply Theorem~\ref{theorem.tightness}, to
conclude that $(\mu_\cdot ^k)_{k\in\N}$ is relatively compact in $\D([0,T],(\MM_F(\XX),\textrm{vague topology}))$. Denote
by $\mu_\cdot$ the limit of any convergent subsequence  $(\mu_\cdot
^{k'})_{k'\in\N}$ and notice that  $\mu_\cdot$ is a.s. strongly
continuous. Indeed, we can easily check that
\[
	\sup_{0\leq t\leq T} 
	\sup_{\substack{f \in \CC^2_b(\XX)\\ \Norm{f}_\infty\leq 1}} 
	\Big|
		\crochet{\mu_t ^{k'},f} 
		- 
		\crochet{\mu_{t^-}^{k'},f} 
	\Big|
	\leq 
	\frac{1}{k'}
\]
holds true, which in turn, implies (\ref{eq.continuity.criteria}).

\subsubsection*{Caracterization}

We expect the martingale and residual terms to vanish when passing to
the limit in (\ref{eq.mu.k}). This would yield  the characterization (\ref{eq.xi.caract}) provided some continuity property holds.

\begin{lemma}
\label{eq.lemma.Psi.continuous}
Let $\mu$ taking values in $\CC([0,T],\MM_F(\XX))$ be the limit of any convergent subsequence $(\mu^{k'})_{k'\in\N}$. Then for any fixed $t\leq T$, $f \in \CC^2_b(\XX)$, the function:
\[
  \Psi_t(\zeta) 
  \eqdef  
  \crochet{\zeta_t,f}
  - 
  \crochet{\zeta_0,f} 
  -
  \int_0^t \bll f(\zeta_{s})\,\rmd s 
\]
is continuous at point $\mu(\omega)$ on $\D([0,T],\MM_F(\XX))$, for all $\omega$ a.s.
\end{lemma}

\proof
As $\mu$ is a continuous process, from the characterization of the Skorohod metric (see \cite[Proposition 6.5 Ch. 3]{ethier1986a}) we get: $\crochet{\mu^{k'}_t,f}$ tends to $\crochet{\mu_t,f}$ for all $t$ and $f$ bounded/continuous. Then:
\begin{align*}
  |\Psi_t(\mu^{k'})-\Psi_t(\mu)|
  &\leq
  |\crochet{\mu^{k'}_t-\mu_t,f}|
  +
  |\crochet{\mu^{k'}_0-\mu_0,f}|
  +
  \sum_{\vartriangle=\textrm{\rm d,b,c,g}}
  \int_0^t |\blltria f(\mu^{k'}_s)-\blltria f(\mu_s)|\,\rmd s
\end{align*}
Both ``d'' and ``b'' terms are of the form $\crochet{\tilde f,\mu^{k'}_s-\mu_s}$ which converges to 0. Consider the ``c'' term, we have $\bllc f(\mu^{k'}_s) \to \bllc f(\mu_s)$ because the function $\zeta_s \mapsto \int\int u(x,y)\,f(x)\, \zeta_s\otimes \zeta_s(\rmd x,\rmd y)$ is continuous (for the weak topology).

For the ``g'' term:
\begin{align*}
  &|\bllg f(\mu^{k'}_s)-\bllg f(\mu_s)|
  \leq
  \Big|
    \int_\XX [\bar\GG f(x,\mu^{k'}_s)-\bar\GG f(x,\mu_s)]\,\mu^{k'}_s(\rmd x)
  \Big|
\\
  &\qquad\qquad\qquad\qquad\qquad\qquad\qquad
  +
  \Big|
    \int_\XX \bar\GG f(x,\mu_s) \,[\mu^{k'}_s(\rmd x)-\mu_s(\rmd x)]
  \Big|
\end{align*}
with $\bar\GG(x,\zeta)=\nabla f(x)\cdot\bar g(x,\zeta) + \demi\, 
\frac{\partial^2 f(x)}{\partial x_\ell\,\partial x_{\ell'}}
             \,\bar a_{\ell\ell'}^k (x,\zeta)
             $.
Thanks to the Lipschitz continuity (\ref{eq.Lipschitz}), the first term of the {\sc rsh} of the previous inequality is bounded by:          
\[
  \crochet{\mu^{k'}_s,1} \times \sup_{x\in\XX}  
  |\bar\GG f(x,\mu^{k'}_s)-\bar\GG f(x,\mu_s)|
  \leq
  C\,\crochet{\mu^{k'}_s,1} \, \rho(\mu^{k'}_s,\mu_s)
\]
which tends to 0 as $\mu^{k'}_s\xrightarrow[]{w}\mu_s$. The first term also tends to 0.
\carre

\bigskip

Note that
\[
  \Psi_t(\mu ^k) 
  = 
  Z^k_t + R^k_t\,.
\]
where $Z^k_t$ and $R^k_t$ are given by (\ref{eq.Z.k}) and (\ref{eq.R.k}) respectively.
$Z^k_t$ is a centered martingale, we now check that its quadratic variation tends to 0. Indeed,
Using (\ref{eq.crochet.Z.k}) and (\ref{eq.lemma.pop.finie.unif}), we have  
\[
  \E ( |Z^{k'}_t |^2) 
  = 
  \E \crochet{Z^{k'}}_t 
  \leq
  \frac{C_f}{k'}\, 
  \E
     \int _0^t \big(\crochet{\mu^{k'}_s,1 } + \crochet{\mu^{k'}_s,1 } ^2\big)
      \, \rmd s 
  \leq 
  \frac{C_{f,t}}{k'}
\]
which shows that $Z^{k'}_t$ goes to~0 in $L^2$. Similarly, since
\[
  \sum_{\vartriangle=\textrm{\rm c,g}}
  \Big|
  \frac{1}{k} \lltria f(k\, \mu) - \bar  \lltria
	f(\mu)
  \Big|
	\leq C_f \,( \crochet{\mu,1} +  \crochet{\mu,1} ^2) \ ,
\]
we conclude that $R^{k'}_t$ goes to~0 in $L^1$ thanks
to~(\ref{eq.lemma.pop.finie.unif}) and Lebesgue's Theorem. Finally, for
$\zeta \in \D([0,T],\MM_F(\XX))$, 
\[
  |\Psi_t(\zeta)| 
  \leq 
  C_{f,t} \, 
  \sup_{0\leq s\leq t} ( \crochet{\zeta_s,1} +  \crochet{\zeta_s,1} ^2)
\]
so that $( \Psi_t(\mu^{k'}))_{k'}$ is uniformly integrable
by~(\ref{eq.lemma.pop.finie.unif}). Lebesgue's theorem and
a.s. continuity of $\Psi$ at $\mu$ yields
\[
  0 
  = 
  \lim_{k'} \E |\Psi_t(\mu ^{k'})| 
  = 
  \E|\Psi_t(\mu)| 
  \,.
\]
Using asumption (\ref{eq.condi.ini}), we conclude that $\mu^k$ converges to a process which is a.s. equal to the deterministic unique solution of~(\ref{eq.xi.caract}).

\subsubsection*{Weak convergence}

So far, we have established the convergence of  $\mu^k$ to $\xi$ in
$\D([0,T],(\MM_F(\XX),\textrm{vague topology}))$. The extension to the weak topology 
is achieved thought \cite[Th. 3]{meleard1993a}. Indeed, we only need to prove that
the total mass $\crochet{\mu_\cdot ^k,1}$ converges in law to
$\crochet{\xi_\cdot ,1}$ in $\D([0,T],\R)$, provided that the limitting
process is continuous. But this is a particular
case of the work done with $f=1$.

\bigskip

This ends the proof of Theorem \ref{th.cv}.

\section{Conclusion}

Using the approach described in this paper, most IBMs could be rewritten as measure-valued Markov branching/diffusing processes. This allows first an insight into the Monte Carlo procedure adapted to the situation, and second to derive a model at macroscopic level as a limit model on various asymptotic situations.

The model presented in this work assumes that there is no limitation for resources. For concrete applications, this model should be coupled with models of the resource dynamics. It is also important to devote more work to the case of ecosystems in the presence of several species.

From the simulation point of view, the direct microscopic IBM approach is limited to relatively small population sizes. For realistic scenarios it could be possible to use a coupled approach: species of importance with small population sizes, could be modeled with IBMs; while less important species with large population sizes, could be modeled through coarser macroscopic models.

\appendix
\section*{Appendix}
\addcontentsline{toc}{section}{Appendix}

\section{Weak convergence of sequence of measure-valued processes}

Consider the space $\D([0,T],\MM_F(\XX))$ of \emph{càdlàg}  functions with values in the space $\MM_F(\XX)$ of finite positive measures. The space $\D([0,T],\MM_F(\XX))$ is endowed with the classical Skorokhod topology. The space $\MM_F(\XX)$ will be equipped with the weak or the vague topologies: the weak (resp. vague) topology corresponds to the convergence $\nu_n\xrightarrow[]{\textrm{\tiny weak}}\nu$
(resp. $\nu_n\xrightarrow[]{\textrm{\tiny vague}}\nu$) defined by:
\[
   \crochet{\nu_n,\phi}
   \xrightarrow[n\to\infty]{}
   \crochet{\nu,\phi}
   \qquad
   \forall \phi\in\CC_b(\XX)
   \  \textrm{(resp. $\CC_K(\XX)$)}
\]
where $\CC_b(\XX)$ is the set of bounded continuous functions from $\XX$ to $\RR$ (resp. $\CC_K(\XX)$ is set of continuous functions from $\XX$ to $\R$ with compact support). Hence the Skorokhod topology on $\D([0,T],\MM_F(\XX))$ could be associated  to the weak or to the vague topology on $\MM_F(\XX)$. Unless otherwise specified, $\D([0,T],\MM_F(\XX))$ will be equipped with the weak topology.

\bigskip

Consider a sequence of random variables $(\mu^k)_{k\in \N}$ which take values in $\D([0,T],\MM_F(\XX))$ (i.e. $(\mu^k_t)_{t\in[0,T]}$ is a \emph{càdlàg} process which takes values in $\MM_F(\XX)$). Here we want to summarize the conditions that ensure that $\mu^k$ converges in law, that is $Q^k$ converges weakly where $Q^k$ denotes the law of $\mu^k$ on $\D([0,T],\MM_F(\XX))$. 
If $Q$ is the limit law and $\mu$ is a process with law $Q$, then
\begin{align}
\label{eq.weak.cv}
  &
  Q^k \xrightarrow[k\to \infty]{\textrm{weak}}Q
\end{align}
means $\E\Phi(\mu^k) \to \E\Phi(\mu)$ for all function $\Phi:\D([0,T],\MM_F(\XX))\mapsto\R$ continuous and bounded. Hence (\ref{eq.weak.cv}) is:
\begin{align*}
  &
  \mu^k \xrightarrow[k\to \infty]{\textrm{law}}\mu
\end{align*}
This property is sometimes called weak convergence of the processes $\mu^k$ toward $\mu$ \cite{joffe1986a}.

\bigskip

\newcommand{\si} {\textrm{(s$_1$)}}
\newcommand{\sii}{\textrm{(s$_2$)}}

The study of the convergence of $Q^k$ is usually accomplished in two steps \cite{joffe1986a}:
\begin{enumerate}

\item[\si] \emph{The tightness step:} First one shows that the sequence $Q^k$ is relatively compact\footnote{i.e. with compact closure, it is equivalent to the fact that from any subsequence one can extract a convergent subsequence.}. It is equivalent to the uniform tightness property: for every $\epsilon>0$ there exists a compact set $\KK_\epsilon$ in $\D([0,T],\MM_F(\XX))$  such that $\inf_k Q^k(\KK_\epsilon) \geq 1-\epsilon$ (\footnote{According to Prohorov's theorem: Tightness implies uniform relative compactness, and it is equivalent because $\D([0,T],\MM_F(\XX))$ with the Skorohod metric is complete and separable (see \cite[Theorem 2.2 Ch. 3]{ethier1986a}).}). By extension, we will say that the sequence $\mu^k$ is tight if the sequence of their laws $Q^k$ is tight.

\item[\sii] \emph{The limit uniqueness step:}
One proves that there is a set of properties that are fulfilled by the limit of any convergent subsequence $Q^{k'}$ of $Q^{k}$. By showing that only one law $Q$ satisfies this set of properties, we both prove that the sequence $Q^{k}$ converges and characterize the limit law $Q$ with the associated limit process $\mu$. This set of properties is usually expressed in terms of a martingale problem.
\end{enumerate}

\subsection*{The tightness step \textrm{(s$_1$)}}

The tightness step \si\ is achieved through the: 
\begin{theorem}[{\cite[Theorem 9.1 Ch. 3]{ethier1986a}}]
\label{theorem.tightness}
Suppose the the following compact containment condition holds: for every $\epsilon >0$ there exists a compact set $K_\epsilon$ of $\MM_F(\XX)$ such that
\begin{align}
\label{eq.compact.containment.condition}
   \inf_k \P\Big(
      \mu^k_t\in K_\epsilon\ \textrm{for }0\leq t\leq T
   \Big)
   \geq 1-\epsilon\,.
\end{align}
Let $\Theta$ be a dense subset of $\CC_b(\MM_F(\XX);\R)$ for the topology of uniform convergence on compact sets. Then $(\mu^k_.)_{k\in\N}$ is relatively compact if an only if $(f(\mu^k_.))_{k\in\N}$ is relatively compact as a family of processes in  $\D([0,T],\R)$ for any $f\in\Theta$.
\end{theorem}

An example of such a set $\Theta$ is given in \cite[Th. 2.1]{roelly1986a} by $\Theta=\{f_\phi\,;\, \phi\in\tilde\Theta\}$, with $f_\phi(\nu)=\crochet{\nu,\phi}$, where $\tilde\Theta$ is any set dense in $\CC_b(\XX)$.

\bigskip

For any given $f\in\Theta$, to check that:
\[
    Y^k_. \eqdef f(\mu^k_.)
\]
is tight in  $\D([0,T],\R)$. To this end, when $(Y^k_t)_{t\in[0,T]}$ is a semimartingale, we make use of the Aldous-Rebolledo criteria that can be found in \cite[Cor. 2.3.3]{joffe1986a}:
\begin{theorem}[Aldous-Rebolledo criteria]
Let $(Y^k_.)_{k\in \N}$ be a sequence of real valued semimartingales with \emph{càdlàg} paths. If:
\begin{enumerate}
\item 
For each fixed $t\in[0,T]$,  $(Y^k_t)_{k\in \N}$ is tight.

\item 
Let $Y^k_t=A^k_t+M^k_t$ where $A^k_t$ is a finite variation process and $M^k_t$ is a locally square-integrable martingale. Suppose that for any given sequence of stopping times $\tau_k$, bounded by $T$, for each $\epsilon>0$ there exists $\delta>0$ and $k_0$:
\begin{align}
\label{eq.AR.2.1}
  \sup_{k\geq k_0} \sup_{\theta\in[0,\delta]}
  \P\Big(
    \big|
       A^k_{\tau_k+ \theta} - A^k_{\tau_k}
    \big|
    > \epsilon
  \Big)
  &
  \leq \epsilon\,,
\\
\label{eq.AR.2.2}
  \sup_{k\geq k_0} \sup_{\theta\in[0,\delta]}
  \P\Big(
    \big|
       \crochet{M^k}_{\tau_k+ \theta} - \crochet{M^k}_{\tau_k}
    \big|
    > \epsilon
  \Big)
  &
  \leq \epsilon\,.
\end{align} 
\end{enumerate}
Then the sequence $(Y^k_.)_{k\in \N}$  is tight.
\end{theorem}

\bigskip

In order to get a compact set  $K_\epsilon$ of $\MM_F(\XX)$ for the  compact containment condition (\ref{eq.compact.containment.condition}), we can consider a set of the form:
\begin{align}
\label{eq.compact.containment.condition.K}
  K_\epsilon
  &\eqdef
  \Big\{
     \mu\in\MM_F(\XX) \,;\, \crochet{\mu,1} \leq C_\epsilon
  \Big\}\,.
\end{align}
When $\XX$ is not compact, $K_\epsilon$ is compact for the vague topology \cite[Th. 13.4.2]{dieudonne1968a} and not compact for the weak topology. When $\XX$ is compact, the vague topology and the weak topology are the same, and $K_\epsilon$ is compact.

\medskip

\begin{enumerate}

\item
A first approach proposed in \cite{etheridge2000a} is to consider the compactification $\bar\XX$ of $\XX$. Then $\bar K_\epsilon$, the set (\ref{eq.compact.containment.condition.K}) where $\XX$ is replaced by $\bar \XX$, is weakly compact and we get the tightness on $\D([0,T],(\MM_F(\bar\XX),\textrm{weak topology}))$. To conclude, we need to prove that the limit process lies on $\MM_F(\XX)$, i.e. $Q$ is supported by $\D([0,T],\MM_F(\XX))$.

\medskip

\item
A second approach proposed in \cite{roelly1986a} is to work with the vague topology: we first prove tighness on $\D([0,T],(\MM_F(\XX),\textrm{vague topology}))$
and a the associated convergence result. To conclude, according to \cite[Th. 3]{meleard1993a}, if:
\begin{align}
\label{eq.RC93a}
  &
  \mu^k_. \xrightarrow[k\to \infty]{\textrm{law}}\mu_.
  \quad
  \textrm{in }
  \D([0,T],(\MM_F(\XX),\textrm{vague topology}))
  \,,
\\
\label{eq.RC93b}
  &
  \crochet{\mu^k_.,1} \xrightarrow[k\to \infty]{\textrm{law}}\crochet{\mu_.,1}
  \quad
  \textrm{in }
  \D([0,T];\R)
  \,,
\\
\label{eq.RC93c}
  &
  \mu_. \textrm{ is a process on } \CC([0,T];\MM_F(\XX))\,,
\end{align}
then (\ref{eq.weak.cv}) holds true.

According to \cite[Prop. 3.26]{jacod1987a}, to ensure (\ref{eq.RC93c}) we prove that for all $\epsilon>0$:
\begin{align}
\label{eq.continuity.criteria}
  \lim_{k\to\infty}\P\Big(
      \sup_{t\in[0,T]} 
      \big|
        \crochet{\mu^k_t,f}-\crochet{\mu^k_{t^-},f}
      \big|
      \geq \epsilon
  \Big)
  =0
\end{align}
for all $f\in\CC_b(\XX)$.

\end{enumerate}

\subsection*{The limit uniqueness step \sii}

For the step \sii: let $Q$ be the limit of a any convergent subsequence $Q^{k'}$ and $\mu_.$ a process with law $Q$. If we prove that $\mu_.$ satisfies an integro-differential equation that admits a unique solution, then  $Q^{k}\to Q$.

Usually, $Q$ is characterized as a solution of a martingale problem which appears to be the limit of martingale problems satisfied by $Q^k$. At this level, an argument allowing us to pass to the limit in the martingale problem should be invoked as the one proposed in \cite[Ch. 4 Th. 8.10]{ethier1986a}. Then the uniqueness of the solution of the limit martingale problem could be obtained by a duality argument as proposed in \cite[Section 1.6]{etheridge2000a} or \cite[p. 188]{ethier1986a}.

\bibliographystyle{plain}
\bibliography{paper}

\end{document}